\documentclass[journal,onecolumn]{IEEEtran}

\usepackage{cite}
\usepackage{graphics}
\usepackage{amsmath,amssymb,amsfonts}
\usepackage[english]{babel}
\usepackage{graphicx}
\usepackage{color}

\newcommand{\R}{\mathbb{R}}

\newcommand{\tr}{^\top}
\newcommand{\sgn}{{\rm sign}}

\newcommand{\pui}[2]{\left\{{#1}\right\}^{\left\{{#2}\right\}}}

\numberwithin{equation}{section}

\usepackage{color}
\definecolor{Wilfrid}{rgb}    {0.7,0.6,0.0}
\definecolor{Jean-Michel}{rgb}   {0.9,0.7,0.0}
\definecolor{Brigitte}{rgb}   {0.0,0.7,0.6}

\begin{document}
\title{Small-time stabilization of nonholonomic or underactuated mechanical systems: the unicycle and the slider examples}
\author{Brigitte d'Andr\'{e}a-Novel, Jean-Michel Coron, Wilfrid Perruquetti
			\thanks{``This work was supported in part by ANR Finite4SoS (ANR 15 CE23 0007).'' }
			\thanks{Brigitte d'Andr\'{e}a-Novel is with MINES ParisTech-CAOR, PSL Research University, 60 Bvd St-Michel, 75006 Paris, {\tt\small brigitte.dandrea-novel@mines-paristech.fr}}
			\thanks{Jean-Michel Coron is with Sorbonne Universit\'{e}, Universit\'{e} Paris-Diderot SPC, CNRS, INRIA, Laboratoire Jacques-Louis Lions, LJLL,  \'{e}quipe CAGE, F-75005 Paris, France, {\tt\small coron@ann.jussieu.fr}}
			\thanks{Wilfrid Perruquetti is with 
				CRIStAL (UMR-CNRS 9189), Ecole Centrale de Lille, Cit\'ee Scientifique, 59651 Villeneuve-d'Ascq, France, {\tt\small wilfrid.perruquetti@ec-lille.fr}}
			}
		
		\maketitle
		
		\begin{abstract}
			This paper concerns the small-time stabilization of some classes of mechanical systems which are not stabilizable by means of at least continuous state feedback laws. This is the case of nonholonomic mechanical systems, an example being the unicycle robot, or for underactuated mechanical systems, an example being the slider. Explicit time-varying feedback laws leading to small-time stabilization are constructed for these two control systems. The main tools are homogeneity, backstepping, and desingularization technics.
		\end{abstract}
		

\newtheorem{remark}{Remark}[section]
\newtheorem{definition}{Definition}[section]
\newtheorem{theorem}{Theorem}[section]
\newtheorem{lemma}{Lemma}[section]
\newtheorem{proposition}{Proposition}[section]
\newtheorem{corollary}{Corollary}[section]

\section{Introduction}\label{sec-intro}
We consider nonholonomic or underactuated mechanical systems for which asymptotic stabilization cannot be achieved through continuous state feedback laws since they do not satisfy the necessary condition for feedback stabilization due to Brockett (see \cite{1983-Brockett-PM}, \cite[Theorem 11.1]{2007-Coron-AMS}; see also \cite{1990-Coron-SCL} for a slightly stronger necessary condition).

In this paper we address the small-time stabilization problem for both systems, and the approach will be illustrated through two examples: the first one is the so-called ``unicycle'' mobile robot which is a nonholonomic vehicle, and the second one is the ``slider'' which is an under-actuated mechanical system. Our construction of feedback laws stabilizing in small time these two control systems relies on three main ingredients homogeneity, backstepping, and desingularization.

Our paper is organized as follows. In Section \ref{sec-modeling} we briefly recall modeling issues for these two mechanical systems and the similarity between these two kinds of dynamical systems will be emphasized, with respect to controllability and stabilizability. In Section \ref{sec-defs-small-time} we recall some definitions and give results concerning small-time stabilization and homogeneous control systems, and small-time stabilization of the double integrator. We present our feedback laws  stabilizing  in small time the unicycle robot in Section \ref{sec-unicycle} and, in Section \ref{sec-slider}, the slider.

\section{Modeling, controllability and stabilizability properties}\label{sec-modeling}
In this section the dynamical behaviors of  the unicycle robot and the slider are recalled.
\vspace{5mm}
\noindent
\paragraph{The unicycle robot}

\begin{figure}[!htbp]
 \begin{center}
 \includegraphics[width=0.4\textwidth]{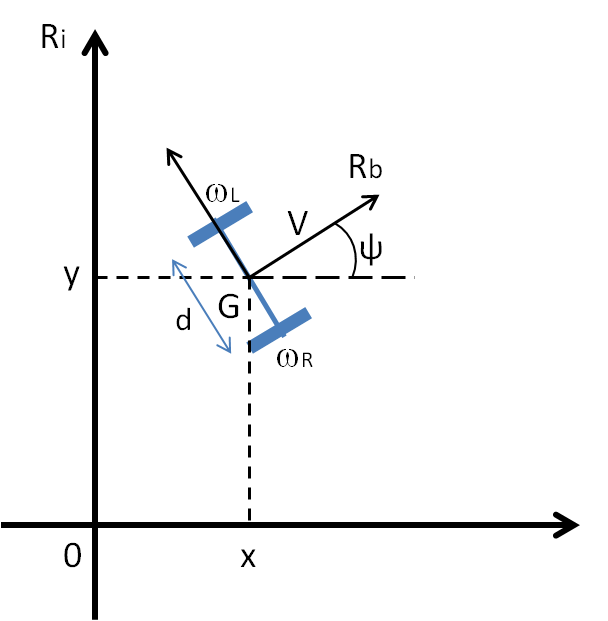}
 \end{center}
 \caption{The unicycle robot}
 \label{fig:unicycle}
\end{figure}

The unicycle is a classical example of non holonomic vehicle. It moves on the 2D-horizontal plane. The configuration vector $q = (x, y, \psi)\tr$ is of dimension $n = 3$. The instantaneous non constrained velocity $\eta = (v_1, \dot{\psi})\tr$ has dimension $2$. It is made of the longitudinal velocity $v_1$ and the angular velocity $\dot{\psi}$. The control variables are the angular rotations $\omega_R$ and $\omega_L$ of the right and left wheels, see Figure (\ref{fig:unicycle}). The third wheel is a free wheel which does not restrict the robot mobility. The nonholonomic property of the unicycle comes from the fact that the robot cannot instantaneously move in the lateral direction. Due to the rolling without slipping assumption, the kinematic behavior of point $G$ can be easily obtained from $\omega_R$, $\omega_L$ (angular velocities of the right and left  wheels), the wheels' radius $R$, and $d$:

\begin{equation}
\left\{
\begin{array}{ccc}
v_1 & = & R (\omega_R + \omega_L),\\
\Omega & = & \frac{R}{2d} (\omega_R - \omega_L).
\end{array}
\right.
\end{equation}
The kinematic equations of the unicycle robot can then be written as follows.

\begin{equation}
\left\{
\begin{array}{lcl}
\dot{x} & = & \cos(\psi)v_1, \\
\dot{y} & = & \sin(\psi) v_1, \\
\dot{\psi} & = & \Omega.
\end{array}
\right.
\label{eq:unicycle}
\end{equation}

\vspace{5mm}
\noindent
\paragraph{The slider}

The slider is an underactuated vehicle, similar to marine vehicles (such as hovercrafts or surface vessels), or to terrestrial quadrotors. For more informations about general underactuated vehicles, see, in particular, \cite{2016-d-Andrea-Novel-Thorel-COCV} and the references therein. 

\begin{figure}[!h]
	\begin{center}
		\includegraphics[width=0.4\textwidth]{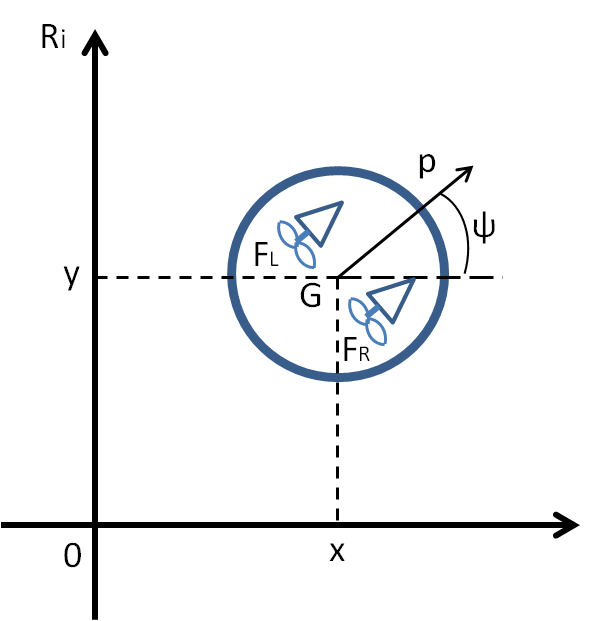}
	\end{center}
	\caption{The slider}
	\label{fig:glisseur}
\end{figure}

The slider moves on a 2D horizontal plane and has a configuration vector $q$ similar to the unicyle's one: $q = (x, y, \psi)\tr$ and $n = 3$. But, the instantaneous velocities of vector $\eta$ are not restricted and $\eta$ has dimension $3$: $\eta = (v_1, v_2, \dot{\psi})\tr$, where $v_1$ and $v_2$ are respectively the longitudinal and lateral velocities in the robot frame, and $\dot \psi$ the angular velocity. Therefore, the control vector has dimension $p = 2 \leq n = 3$. This vehicle is described in Figure (\ref{fig:glisseur}). It is actuated by two propellers producing forces $F_L$ and $F_R$. The sum of these two forces is directly linked to the acceleration of the vehicle, whereas the difference acts on the angular dynamics. Let us denote $\tau_1 = F_L + F_R$ and $\tau_2 = F_R - F_L$, the dynamics can be written:

\begin{equation}
\left\{
\begin{array}{lcl}
\dot{q} & = & S(q)\eta,\\
m \dot{v_1} & = & m v_2 \dot{\psi} + \tau_1,\\
m \dot{v_2} & = & -m v_1 \dot{\psi},\\
I \ddot{\psi} & = & \tau_2,
\end{array}
\right.
\label{eq:sousactionne}
\end{equation}
with
\begin{displaymath}
S(q) =
\left(
\begin{array}{ccc}
\cos (\psi ) & - \sin (\psi ) & 0\\
\sin (\psi ) & \cos (\psi ) & 0\\
0 & 0 & 1
\end{array}
\right).
\end{displaymath}
$m$ being the vehicle mass, and $I$ the vehicle inertia. It should be noticed that friction terms have been neglected. Moreover, $S(q)$ is nothing but the rotation matrix from the inertial frame $R_i$ to the robot frame $R_b$. Therefore, equations (\ref{eq:sousactionne}) describe the dynamic behavior of the slider expressed in the vehicle frame and can be written as follows in the inertial frame:
\begin{equation}
\left\{
\begin{array}{lcl}
m\ddot{x} & = &  \cos (\psi ) \tau_1,\\
m\ddot{y} & = &  \sin (\psi ) \tau_1,\\
I \ddot{\psi} & = & \tau_2.
\end{array}
\right.
\label{eq:sousactionneinertiel}
\end{equation}
The similarity between these equations and the unicycle's ones (\ref{eq:unicycle}) is then obvious.

\subsection{Controllability properties}

The controllability properties of the two systems have been analyzed since a long time. It follows from the Chow-Rashevski theorem (see e.g. \cite[Theorem 3.18]{2007-Coron-AMS}) that the unicycle is globally controllable in finite time as well as small-time locally controllable at any equilibrium of the form $(a,0)\in \R^3\times\R^2$. The global controllability in small time as well as the small-time local controllability of the slider follows from these two controllability properties of the unicycle.  (The small-time local controllability also follows from a general theorem due to Sussmann \cite{1987-Sussmann-SICON}.)
\subsection{Stabilizability properties}
Let us recall that the unicycle control system \eqref{eq:unicycle} and the slider control system \eqref {eq:sousactionneinertiel} cannot be locally asymptotically stabilized  by means of continuous state feedback laws, as it follows  from the following theorem
due to Brockett \cite{1983-Brockett-PM}.

\begin{theorem}
\label{th:Brockett}
Assume that $f:\R^n\times\R^m\to \R^n$ is continuous and satisfies $f(0,0)=0$. Assume that there exists a continuous feedback law $u:\R^n\to\R^m$ vanishing at $0\in \R^n$ such that $0\in \R^n$ is asymptotically stable for the closed-loop system $\dot =f(x,u)$.  Then, for any neighborhood ${\cal N}$ of $(0,0) \in \mathbb{R}^n \times \mathbb{R}^m$
\begin{equation}\label{cond-BBockett}
f({\cal N} ) \text{ is a neighborhood of } 0 \in \mathbb{R}^n.
\end{equation}
\end{theorem}

\begin{remark}\label{Rk-Brockett}
The unicycle robot and the slider do not satisfy the Brockett condition \eqref{cond-BBockett} (and therefore, by Theorem~\ref{th:Brockett}, are not stabilizable by means of at least continuous state feedback laws). Indeed, for the unicycle,  let us consider  $\{ e = (0,\epsilon, 0)^{\top}\}$ with $\epsilon \not = 0$ ; it is clear that this element of $\R^3$ does not belong to the image of $B\times\R^2$ where $B\subset \R^3$ is a ball centered at $0\in \R^3$ of radius strictly  less  $\pi/2$ by the unicycle dynamics given by \eqref{eq:unicycle}. %$\{ f = (U_1, -e_1 \Omega, \Omega)^{\top}\}$ defined at equation \eqref{eq:unicycle_point_fixe}.
In the same way, for the slider, $\{ e = (0, 0, 0, \epsilon, 0, 0)^{\top}\}$ with $\epsilon \not = 0$  does not belong to the image of $ B\times\R^2$ by the slider dynamics given by \eqref{sliderxcossin} below  if  $B\subset \R^6$ is a ball centered at $0\in \R^6$ of radius strictly  less  $\pi/2$.
\end{remark}

Consequently, fixed point asymptotic stabilization needs other control approaches than classical continuous state feedback laws which allow stabilization of exciting reference trajectories. Among those control methods one can cite:
\begin{itemize}
	\item continuous (with respect to the state) time-varying feedback laws,
	\item discontinuous feedback laws.
\end{itemize}

Concerning continuous (with respect to the state) time-varying feedback laws, Coron has shown in \cite{1992-Coron-MCSS} that  asymptotic stabilization of the origin of a small-time locally controllable driftless affine control system can be achieved through smooth time-periodic state feedback laws (see also \cite[Chap. 11, thm. 11.14]{2007-Coron-AMS}). Moreover, this result can be extended to asymptotic stabilization of the origin of many general nonlinear control system  using  continuous time-varying feedback laws; see  \cite{1995-Coron-SICON}, \cite[Chap. 11, thm. 11.28]{2007-Coron-AMS}.
In that context, many authors have proposed stabilizing continuous time-varying feedback laws, see for example \cite{1990-Samson-LN, 1992-Coron-d-Andrea-Nolcos, 1992-Pomet-SCL, 1995-Samson-IEEE, 2000-Morin-Samson-IEEE, 2015-Coron-ICIAM, 2017-Coron-Rivas-SICON, 2016-Guilleron-PhDthesis}  for the unicycle and related control systems (as chained systems). See also \cite{1995-Morin-Samson-Pomet-Jiang-SCL, 1996-Coron-Kerai-A, 1997-Morin-Samson-IEEE} for the stabilization of the attitude of a rigid spacecraft with two controls and \cite{1996-Pettersen-Egeland-CDC, 1997-Pettersen-Egeland-CDC, 2002-Mazenc-Pettersen-Nijmeijer-IEEE, 2017-Coron-Rivas-SICON} for stabilization of underactuated surface vessels. Note however that among these articles the only one dealing with finite-time stabilization is \cite{2016-Guilleron-PhDthesis} .

Concerning discontinuous feedback laws, let us mention \cite{1997-Clarke-Ledyaev-Sontag-Subbotin-IEEE, 2002-Ancona-Bressan-SICON} for  general results showing that controllable systems can be stabilized by means of discontinuous feedback laws. For robustness issues of discontinuous feedback laws,  let us refer to \cite{1999-Sontag-COCV, 2000-Clarke-Ledyaev-Rifford-Stern-SICON, 2004-Ancona-Bressan-SICON, 2005-Prieur-SICON}.  For example of
discontinuous feedback laws for the unicycle or related control systems (as chained systems), let us mention \cite{2003-Prieur-Astolfi-IEEE, 2005-Prieur-Trelat-MCSS}.

We are now interesting in this paper, to design feedback laws ensuring small-time stabilization for both systems. Before that, let us briefly recall some results concerning small-time stabilization.

\section{Small-time stabilization}
\label{sec-defs-small-time}
Let us first introduce some definitions about stabilization and homogeneity.

\begin{definition}
\label{defSTS}
Let $F\in C^0(\R^n,\R^n)$ be such that $F(0)=0$. The origin $0\in \R^n$ is said to be small-time stable  for $\dot x=F(x)$  if $\forall \varepsilon>0,\; \exists \eta>0$ such that
\begin{gather}
\label{0entempspetit}
\left(\dot x =F(x) \text{ and }  |x(0)| < \eta \right) \Rightarrow \left(x(\varepsilon)=0\right),
\\
\label{property-stability}
\left(\dot x =F(x) \text{ and }  |x(0)|<\eta \right) \Rightarrow \left(|x(t)|<\varepsilon%\WilfridCor{delete =0}
, \; \forall t \in \left[0,+\infty\right) \right).
\end{gather}
\end{definition}
\begin{remark}
\label{rem-autre-stability}
Let us point out that, as one can easily check, \eqref{0entempspetit} and \eqref{property-stability} are equivalent to \eqref{0entempspetit} and
\begin{gather}
\label{reste-a-0}
 \left(\dot x =F(x) \text{ and } x(0)=0 \right) \Rightarrow \left(x(t)=0, \; \forall t \in \left[0,+\infty \right)\right).
\end{gather}
Let us mention that "Small-time stability" is in between the classical concepts of Finite-time stability (see \cite{Moulay06,Polyakov15}) and Fixed-time stability (see \cite{Polyakov12,Polyakov15}). "Small-time stability" is equivalent to finite-time stability when the settling time function is continuous.
\end{remark}

Let $r = (r_1, \cdots , r_n )$ be a $n$-uplet of positive real numbers, thereafter called  a {\it weight}. Then,  let us define $\Lambda_r : (0,+\infty) \times \R^n\to \R^n$ by
\begin{equation}\label{deflambdar}
  \Lambda_r (\lambda, x) = (\lambda ^{r_1} x_1,\; \ldots , \lambda ^{r_i} x_i , \ldots, \lambda ^{r_n} x_n ),
\end{equation}
$\forall \lambda \in (0,+\infty),\; \forall x=(x_1,\ldots, x_i,\ldots, x_n)\tr \in \R^n.$ This map %\WilfridCor{delete "is called"}
$\Lambda_r$ is called a dilation. Following \cite{1991-Hermes-part-book} this map allows to define the notion of functions and vector fields which are  $r$-homogeneous .

\begin{definition}
\label{defhomo1}
Let $\kappa \in \mathbb{R}$. A function $h : \mathbb{R}^n \rightarrow \mathbb{R}$ is said to be $r$-homogeneous of degree $\kappa $ if
\begin{gather}
\label{def_homoge-function}
h (\Lambda_r (\lambda,x)) = \lambda ^{\kappa} h(x),\; \forall \lambda \in (0,+\infty),\;
  \forall x  \in \R^n.
\end{gather}
\end{definition}

\begin{definition}
Let $\kappa \in \mathbb{R}$. A vector field $F=(F_1,\ldots,F_i,\ldots, F_n) : \mathbb{R}^n \rightarrow \mathbb{R}^n$ is said to be $r$-homogeneous of degree $\kappa \in \mathbb{R}$ if, for every $i\in \{1,\dots,n\}$, the coordinate function $F_i$ is $r$-homogeneous of degree $\kappa + r_i$.
\label{defhomo2}
\end{definition}

Among many properties of homogeneous systems, let us recall the following result which will be useful in the sequel.

\begin{theorem}[\cite{1992-Rosier-SCL}]\label{thm2}
Let $F$ be a continuous $r$-homogeneous vector-field on $\mathbb{R}^n$ of degree $\kappa \in \mathbb{R}$. If $0\in \R^n$ is asymptotically stable for $\dot x=F(x)$, then, for every real number $\alpha > \max \{r_1, \ldots, r_n\}$, there exists a function $V : \mathbb{R}^n \rightarrow \mathbb{R}$ of class $C^1$  such that
\begin{gather}
\label{V-homogeneous-thm2}
V \text{ is $r$-homogeneous with degree $\alpha$},
\\
V(x)>V(0)=0, \; \forall x \in \R^n\setminus\{0\},
\\
\label{eq:Vhom-thm2}
\displaystyle \frac{\partial V}{\partial x}\cdot  F(x) <0, \; \forall x \in \R^n\setminus\{0\}.
\end{gather}
\end{theorem}
\begin{corollary}[\cite{1992-Rosier-SCL}]\label{cor2}
Let $F$ be a continuous $r$-homogeneous vector-field on $\mathbb{R}^n$ of degree $\kappa \in (-\infty,0)$. If $0\in \R^n$ is asymptotically stable for $\dot x=F(x)$, then $0\in\R^n$ is small-time stable for $\dot x =F(x)$.
\end{corollary}

Let us now consider the control system
\begin{equation}\label{defSigma}
\dot x = f(x,u)
\end{equation}
where the state is $x\in \R^n$, the control is $u\in \R^m$, and $f$ is continuous and satisfies $f(0,0)=0$.
In the sequel, we intend to design time-piecewise continuous periodic feedback laws for \eqref{defSigma}. Such feedback laws are defined as follows.
\begin{definition}
\label{defTPC}
Let $u : \mathbb{R} \times \mathbb{R}^n \rightarrow \mathbb{R}^l$, $(t,x) \mapsto u(t,x)$ is a
time-piecewise continuous function if there exists a sequence of real numbers $t_i \in \mathbb{Z}$ such that:
\begin{gather}
t_i < t_{i+1} \; , \; \forall \; i \; \in \mathbb{Z},\\
\lim_{i \rightarrow + \infty} t_i = + \infty,\\
\lim_{i \rightarrow - \infty} t_i = - \infty,\\
\begin{array}{c}
u_{\mid (t_i , t_{i+1})\times \mathbb{R}^n} \; \text{is the restriction to}\;(t_i , t_{i+1})\times \R^n \text{of a continuous function on}\;[t_i, t_{i+1}] \times \mathbb{R}^n.
\end{array}
\end{gather}
Moreover these feedback laws are called time-piecewise continuous stationary periodic feedback laws if for $t\in (t_i,t_{i+1})$, $u(t,x) $ does not depend on $t$.
%OK pour \modifsjm{Moreover these feedback laws are called time-piecewise continuous stationary periodic feedback laws if for $t\in (t_i,t_{i+1})$, $u(t,x) $ does not depend on $t$.}
\end{definition}
From now on, we denote by $C_{\text{t-p}} (\mathbb{R} \times \mathbb{R}^n; \mathbb{R}^l)$ the set of time-piecewise continuous functions $u(t,x)$ from $\mathbb{R} \times \mathbb{R}^n$ to $\mathbb{R}^l$.

Let us recall the definition of local asymptotic stability.
\begin{definition}
Let $F \in C_{\text{t-p}} (\mathbb{R} \times \mathbb{R}^n; \mathbb{R}^n)$ be such that
\begin{gather}
F(t,0) = 0 , \; \forall \; t \in \mathbb{R}.
\end{gather}
One says that $0$ is locally asymptotically stable for $\dot x=F(t,x)$ if
\begin{itemize}
\item[(i)] For every $\varepsilon>0$ there exists $\eta>0$ such that, for every $s\in \mathbb{R}$ and for every $\tau \geq s$
\begin{gather}
(\dot x = F(t, x), \; |x(s) |< \eta ) \Rightarrow  (|x(\tau ) |< \epsilon),
\end{gather}
\item[(ii)] there exists $\delta > 0$ such that for every $\epsilon > 0$ there exists $M>0$ such that for every $s \in \mathbb{R}$
\begin{gather}
(\dot x = F (t,x) , \; |x(s) |< \delta) \; \Rightarrow \; (|x (\tau ) |< \epsilon , \; \forall \tau \geq s+M ).
\end{gather}
\end{itemize}
\label{def-as-stab}
\end{definition}

Concerning small-time stabilizability by means of time-piecewise continuous periodic feedback laws, we adopt the
 following definition.
\begin{definition}
The control system \eqref{defSigma} is locally stabilizable in small time by means of time-piecewise continuous periodic feedback laws of arbitrary small period if, for every positive real number $T$, there exist $\epsilon > 0$ and $u \in C_{\text{t-p}} (\mathbb{R} \times \mathbb{R}^n; \mathbb{R}^m)$ such that:
\begin{gather}
\label{unulen0}
u(t,0) = 0 , \; \forall \; t \in \mathbb{R},\\
\label{uperiodic}
u(t+T,x) = u(t,x) , \; \forall \; t \in \mathbb{R} \;, \; \forall \; x \in \mathbb{R}^n,\\
\label{0toujours}
(\dot x = f(x,u(t,x)) \;\text{and} \; x(s)=0) \Rightarrow %\nonumber\\
 (x(\tau ) = 0, \; \forall \tau \geq s), \; \forall s \in \mathbb{R},\\
\label{epsilontoujours}
(\dot x = f(x,u(t,x)) \;\text{and} \; |x(s)|\leq \epsilon) \Rightarrow %\nonumber\\
(x(\tau) = 0, \; \forall \tau \geq s+T), \; \forall s \in \mathbb{R}.
\end{gather}
\label{defstabTP}
\end{definition}

\begin{remark}
If the control $u$ is as in Definition~\ref{defstabTP} then the origin of $\mathbb{R}^n$ is locally asymptotically stable for the closed-loop system $\dot x = f(x,u(t,x))$.
\end{remark}

%Ajouter stabilization double integrateur qui sera utile pour le slider..
Before studying the small-time stabilization problem for the unicycle and the slider, let us first prove a result concerning cascaded-systems and let us then apply it first to the small-time stabilization of the double integrator.

\subsection{Small-Time stabilization for homogeneous cascaded systems}
Let us consider the following control system:

\begin{equation}
\dot x = f(x,y), \;\dot y = v,
\label{eqcascaded}
\end{equation}
where the state is $(x,y)\in \R^n\times \R$ and the control is $v\in \R$. Let
$r=(r_1,\ldots, r_n,r_{n+1})\tr \in (0,+\infty)^{n+1}$. For $x=(x_1,\ldots,x_n)\tr \in \R^n$ and for
$\lambda \in (0,+\infty)$, we  use the notation $\Lambda_{\bar r} (\lambda, x)$ introduced in (\ref{deflambdar}) with $\bar r := (r_1, \cdots , r_n )$.
%\begin{equation}\label{defdeltalambda}
%\delta_\lambda(x)=(\lambda^{r_1}x_1,\ldots,\lambda^{r_n}x_n).
%\end{equation}
We assume that
$f=(f_1,f_2,\ldots,f_n)$ is of class $C^0$ on $\R^n\times \R$ and that, for every $i\in \{1,\ldots,n\}$, $f_i$ satisfies the following homogeneity property
\begin{equation}
\label{eq:prophom1}
f_i (\Lambda_{\bar r} (\lambda, x), \lambda^{r_{n+1}} y ) = \lambda^{r_i + \kappa} f_i (x,y), 
\end{equation}
$\forall x\in \R^n, \; \forall y \in \R, \forall \lambda \in (0,+\infty)$. Let us define $\pui{x}{\alpha}={\sgn}(x)|x|^{\alpha}$, for every real number $\alpha > 0$
and for every real number $x$. The next theorem is dealing with backstepping and stabilization by means of a continuous feedback law.
\begin{theorem}\label{thmcascaded}
Let us suppose that $\kappa\in \R$ satisfies
	\begin{equation}
	\kappa  + r_{i}>0, \forall i \in \{1,\ldots,n,n+1\},
	\label{kappagd}
	\end{equation}
and that there exist a continuous feedback law
	$\bar y : \mathbb{R}^n \rightarrow \mathbb{R}$ and $l\in (0,+\infty)$ such that
\begin{gather}
\label{lgrand}
l+1>\displaystyle \frac{r_i}{r_{n+1}}, \; \forall i \in \{1,\ldots,n\},
\\
\label{homog-bary}
\text{$\bar y$ is  $\bar r$-homogeneous of degree $r_{n+1}$},
\\
\label{asforbary}
\text{$0\in \R^n$ is asymptotically stable for $\dot x = f(x, \bar y (x))$,}
\\
\label{new-assumption-l}
\text{$\pui{\bar y}{l}$ is of class $C^1$}.
\end{gather}
Then the control law
\begin{equation}
\label{feedback-v(x,y)}
v (x,y) = - k \pui{\pui{y}{l} - \pui{\bar y (x)}{l}}{(\kappa + r_{n+1})/(lr_{n+1})},
\end{equation}
with a sufficiently high gain $k>0$, makes $0\in \R^{n+1}$  asymptotically stable for the cascaded system \eqref{eqcascaded}. Moreover, if $\kappa < 0$, this asymptotic stability is a small-time stability.
\end{theorem}

\begin{remark}
\label{increasingl}
Note that property \eqref{new-assumption-l} implies that, for every $L\in [l,+\infty)$,
\begin{equation}\label{regpourL>l}
\text{$\pui{\bar y}{L}$ is of class $C^1$}.
\end{equation}
Hence, \eqref{lgrand} is not an important restriction: it can always been achieved by increasing $l$ if necessary.

\end{remark}

\begin{remark}
\label{rem-Morin-Samson}
Theorem~\ref{thmcascaded} is related to  \cite[Proposition 5]{1997-Morin-Samson-IEEE}. The main improvement of
Theorem~\ref{thmcascaded} compared to  \cite[Proposition 5]{1997-Morin-Samson-IEEE} is that
Theorem~\ref{thmcascaded} can be applied in a recursive manner, which is not the case of \cite[Proposition 5]{1997-Morin-Samson-IEEE} since, with the notations of this proposition, the feedback law
$(x_1,y,t)\in \R^m\times \R\times \R \mapsto y-v(x_1,t)\in \R$ is not necessarily of class $C^1$ on
$((\R^m\times \R)\setminus\{(0,0)\})\times \R$. To apply Theorem~\ref{thmcascaded} in a recursive manner, let us first point out that, by Remark~\ref{increasingl}, increasing $l$ if necessary,  we may assume that $l\geq 1$. %OK \modifsjm{by Remark~\ref{increasingl}, increasing $l$ if necessary,  we may assume that $l\geq 1$}. 
Then $v$ defined by \eqref{feedback-v(x,y)} is such that
\begin{gather}
\label{vpowerC1}
\pui{v}{(lr_{n+1})/(\kappa + r_{n+1})} \text{ is of class $C^1$.}
\end{gather}
%OK \commentsjm{C'est pour \eqref{vpowerC1} que l'on a besoin de $l\geq 1$.} 
Let us also point that $v$  satisfies the
following homogeneity property \begin{gather}\label{homogeneity-vxy}
v(\Lambda_{\bar r} (\lambda, x), \lambda^{r_{n+1}} y )=\lambda^{r_{n+1} + \kappa} v (x,y),
\end{gather}
$\forall x\in \R^n, \; \forall y \in \R,
\forall \lambda \in (0,+\infty)$. Let
\begin{equation}\label{defrn+2}
r_{n+2}:= r_{n+1}+\kappa.
\end{equation}
One can then apply Theorem~\ref{thmcascaded} to the control system
\begin{equation}\label{f-extended}
\dot x= f(x,y), \; \dot y = v, \; \dot v=w,
\end{equation}
where the state is $(x,y,v)\in \R^n\times \R\times \R$ and the control is $w\in \R$, provided that (compare to \eqref{kappagd})
\begin{gather}\label{new-condition-l}
r_{n+2}+\kappa=r_{n+1}+2\kappa  >0,
\\
\label{condition-l-cran-d-apres}
l+1>\frac{r_i}{r_{n+1}+\kappa},\; \forall i \in \{1,\ldots,n,n+1\}.
\end{gather}
 One can keep going $k$ times (the state being then in $\R^{n+k+1}$) provided that
\begin{equation}\label{new-condition-k}
r_{n+1}+(k+1)\kappa  >0.
\end{equation}
In particular, if $\kappa\geq 0$, one can keep going as long as desired. However in the case where
small-time stability is desired (i.e. the case where $\kappa<0$), condition \eqref{new-condition-k} provides an upper bound on $k$.
As mentioned in Remark~\ref{increasingl}, the condition on $l$
(see, e.g., \eqref{condition-l-cran-d-apres}) can always be achieved by increasing $l$.
\end{remark}

 \begin{remark}
 It has been already proved in \cite{1991-Coron-Praly-SCL} that under the assumptions of Theorem~\ref{thmcascaded} and even without the assumption \eqref{new-assumption-l}, there exists a continuous feedback law satisfying the homogeneity
 property \eqref{homogeneity-vxy}. Compared to \cite{1991-Coron-Praly-SCL} the interest of Theorem~\ref{thmcascaded} is that it provides
 an explicit asymptotic stabilizing feedback law. For $l=1$ Theorem~\ref{thmcascaded} is already known: see \cite[Prop. 5]{1997-Morin-Samson-IEEE}.
\end{remark}
\noindent
\textbf{Proof of Theorem~\ref{thmcascaded}.} Let $\bar \rho : \R^n \to [0,+\infty)$ be the following homogeneous ``norm'' associated to the dilation $(x_1,\ldots,x_n)\tr \in \R^n \to
(\lambda^{r_1}x_1,\ldots,\lambda^{r_n}x_n)\tr \in \R^n$:
\begin{equation}
\bar \rho (x) = \sum_{i=1}^n |x_i |^{1/r_i}, \; \forall x=(x_1,\ldots,x_n)\tr \in \R^n.
\label{defrhobar}
\end{equation}
Let $\alpha \in (0,+\infty)$ be defined by
\begin{equation}\label{defalpha}
\alpha:=(l+1)r_{n+1}.
\end{equation}
 Note that, from \eqref{lgrand}, we have
\begin{equation}\label{alpha-assez-grand}
\alpha>r_i, \; \forall i \in \{1,\ldots,n\}.
\end{equation}
Using Theorem~\ref{thm2} with $r:=\bar r$ and $F(x):=f(x,\bar y (x))$, the assumptions of Theorem~\ref{thmcascaded} and standard homogeneity arguments, one gets the existence of  $V\in C^1(\R^n,\R)$ and $c_1\in (0,1)$ such that
\begin{gather}
\label{V-homogeneous}
V \text{ is $\bar r$-homogeneous with degree $\alpha>0$,}
\\
c_1 (\bar \rho (x))^{\alpha} \leq V(x)\leq \frac{1}{c_1} (\bar \rho (x))^{\alpha}, \; \forall x \in \R^n,
\\
\label{eq:Vhom}
\displaystyle \frac{\partial V}{\partial x}(x)\cdot  f(x, \bar y (x)) \leq - 2 c_1 (\bar \rho (x))^{\alpha+ \kappa},\; \forall x\in \R^n.
\end{gather}
 Using the ``desingularization method'' proposed in \cite{1991-Praly-dAndrea-Novel-Coron-IEEE}, let us now introduce the following control Lyapunov function candidate for the control system \eqref{eqcascaded}:

\begin{equation}
W(x,y) = V(x) + a  \Phi (x,y),
\label{eq:Whom}
\end{equation}
with $a$ chosen in $(0,+\infty)$ and
\begin{eqnarray}
\Phi (x,y) := \int _{\bar y (x)} ^y \left(\pui{s}{l} - \pui{\bar y(x)}{l}\right) ds 
= \displaystyle \frac{l}{l+1} |\bar y (x)|^{l+1} - y \pui{\bar y (x)}{l} + \displaystyle \frac{|y|^{l+1}}{l+1}.
\label{eq:Whom2}
\end{eqnarray}
Let us compute the time-derivative $\dot W$ along the trajectories of the closed-loop system
\begin{equation}
\dot x = f(x,y), \;
\dot y = - k \pui{\pui{y}{l} - \pui{\bar y (x)}{l}}{(\kappa + r_{n+1})/(lr_{n+1})}.
\label{closed}
\end{equation}
We get
\begin{equation}
\dot W = \dot V + a \dot \Phi =Q_1 + Q_2 + ak Q_3 + a Q_4
\label{Wpoint}
\end{equation}
with
\begin{gather}
Q_1(x) := \displaystyle \frac{\partial V}{\partial x} f(x, \bar y (x)),\\
Q_2 (x,y):= \displaystyle \frac{\partial V}{\partial x} \left ( f(x,y) - f(x, \bar y (x)) \right ),\\
\label{defQ3}
Q_3(x,y) := -\displaystyle \frac{\partial \Phi}{\partial y} \pui{\pui{y}{l} - \pui{\bar y (x)}{l}}{(\kappa + r_{n+1})/(lr_{n+1})},\\
Q_4 (x,y):=  \displaystyle \frac{\partial \Phi}{\partial x} f(x,y).
\end{gather}

Let us point out that $Q_1$, $Q_2$, and $Q_3$  are well defined and continuous.
Noticing that
\begin{equation}\label{expression-abbsyl+1}
|\bar y (x)|^{l+1}= |\pui{\bar y(x)}{l}|^{(l+1)/l},
\end{equation}
 one sees, using also \eqref{new-assumption-l}, that $Q_4$ is also well defined and continuous on $\R^n$.
By \eqref{eq:Vhom}, we have
\begin{equation}
Q_1 \leq - 2 c_1 \bar \rho^{\alpha + \kappa}.
\label{estQ1}
\end{equation}

Let us now deal with $Q_2$. For that purpose, let us define, for $p \in \mathbb{N}^{\star},  (x,y) \in \mathbb{R}^n \times \mathbb{R} \setminus \{(0,0)\}$
\begin{equation}
G_p(x,y) = \displaystyle \frac{\displaystyle \frac{\partial V}{\partial x}(x) \left ( f(x,y) - f(x, \bar y (x)) \right )}{c_1 (\bar \rho (x))^{\alpha + \kappa} + p |\pui{y}{l} -
\pui{\bar y(x)}{l} |^{(\alpha + \kappa)/(lr_{n+1})}}.
\label{defGp}
\end{equation}
Let us prove that if $p$ is large enough, we have
\begin{equation}
|G_p |\leq 1, \;  (x,y) \in \mathbb{R}^n \times \mathbb{R} \backslash \{(0,0)\}.
\label{ineqGp}
\end{equation}
We argue by contradiction, and therefore assume the existence of a sequence
$(x^p, y^p)_{p\in \mathbb{N}} \in \mathbb{R}^n \times \mathbb{R} \setminus \{(0,0)\}$ such that, for $p$ large enough,
\begin{eqnarray}
c_1 (\bar \rho (x^p))^{\alpha + \kappa} +  p |\pui{y^p}{l} -
\pui{\bar y(x^p)}{l} |^{(\alpha + \kappa)/(lr_{n+1})} %\nonumber \\
< \displaystyle \left\vert  \frac{\partial V}{\partial x}(x^p) \left ( f(x^p,y^p) - f(x^p, \bar y (x^p)) \right ) \right\vert.
\label{suite}
\end{eqnarray}
Let us point out that $G_p$ is $r$-homogeneous with degree 0. Hence, we may assume that $\rho (x^p,y^p) = 1$ with
\begin{equation}
\label{defrho}
\rho (x,y ) := | \pui{y}{l}- \pui{\bar y (x)}{l}\mid^{1/(lr_{n+1})} + \bar \rho (x).
\end{equation}
In particular, the sequence $(x^p, y^p)_{p\in \mathbb{N}}$ is bounded and, extracting a subsequence if necessary, we may assume that
$x^p \rightarrow x^{\star}$ and $y^p \rightarrow y^{\star}$ as $p \rightarrow + \infty$, with
\begin{equation}
\rho (x^{\star},y^{\star}) = 1.
\label{rhostar1}
\end{equation}
Note that by \eqref{kappagd} and \eqref{defalpha}, $\alpha +\kappa>0$. Hence, using inequality \eqref{suite}, we get that $y^{\star} = \bar y(x^{\star})$ and also that $\bar \rho (x^{\star}) = 0$. Therefore, $x^{\star}= 0$, and $y^{\star} = 0$ which is in contradiction with \eqref{rhostar1}. Hence, we can now choose $p = c_2$ large enough so that we have:
\begin{equation}
|Q_2(x,y) |\leq c_1 (\bar \rho (x))^{\alpha+\kappa} + c_2 | \pui{y}{l}- \pui{\bar y (x)}{l}\mid^{1/(lr_{n+1})},
\label{estQ2}
\end{equation}
and it should be noticed that $c_2$ does not depend on $a,k$.

Concerning  $Q_3$ we just point out that, using \eqref{defalpha}, \eqref{eq:Whom2}, and \eqref{defQ3}, we have
\begin{equation}
Q_3(x,y) = - | \pui{y}{l}- \pui{\bar y (x)}{l}|^{(\alpha + \kappa)/(lr_{n+1})}.
\label{Q33}
\end{equation}

Let us now deal with the last term $Q_4$. From \eqref{eq:Whom2}, $\Phi$ is $r$-homogeneous of degree $\alpha$ and therefore $Q_4$ is $r$-homogeneous of degree $\alpha+\kappa$. Similarly
\begin{equation}
(\bar \rho (x))^{\alpha + \kappa}+ p\vert y - \bar y (x)\vert^{(\alpha + \kappa)/r_{n+1}}, \; p>0,
\label{Q41}
\end{equation}
is $r$-homogeneous of degree $\alpha +\kappa$ and is strictly positive outside the origin.
Proceeding as for $Q_2$, one can show that there exists a sufficiently large positive constant $c_4$ such that:
\begin{equation}
\vert	Q_4 \vert \leq (\bar \rho (x))^{\alpha + \kappa} + c_4| \pui{y}{l}- \pui{\bar y (x)}{l}|^{(\alpha + \kappa)/(lr_{n+1})}.
	\label{Q42}
\end{equation}
%Then, using \eqref{Wpoint}, \eqref{estQ1}, \eqref{estQ2}, \eqref{Q33} and \eqref{Q42} we can write:
%
%\begin{equation}
%\dot W \leq (-K_1 + K_4) (\bar \rho (x))^{\alpha + \kappa} + (- 2^{1-l} \sqrt{K} + K_2 + K_4 ) (y - \bar y (x))^{(\alpha + \kappa)/r_{n+1}}
%\label{Wpointfin}
%\end{equation}
%For $K$ and $K_1$ sufficiently large, it is clear that $\dot W$ is negative and equal to zero if and only if $x=0$ and $y=0$ which allows to conclude that the origin of the closed-loop system \eqref{closed} is asymptotically stable.
%Moreover, if the homogeneity degree $\kappa$ of the closed-loop vector field is strictly negative, this asymptotic stability is a finite-time stability. This concludes the proof.
%
%
%%%%%
Then, using \eqref{Wpoint}, \eqref{estQ1}, \eqref{estQ2}, \eqref{Q33}, and \eqref{Q42}, we get
\begin{eqnarray}
	\dot W &\leq& (-c_1 +a
	) (\bar \rho (x))^{\alpha + \kappa} + (c_2+ac_4-ak) | \pui{y}{l} \nonumber\\
	&&- \pui{\bar y (x)}{l}|^{(\alpha + \kappa)/(lr_{n+1})}.\label{Wpointfin}
\end{eqnarray}
 For $a<c_1$ and $k$ sufficiently large ($k>(c_2+ac_4)/a$), we get that $\dot W$ is non positive and is equal to zero if and only if $x=0$ and $y=0$ which allows to conclude that the origin of $ \R^n\times \R$ is asymptotically stable for the closed-loop system \eqref{closed}. If $\kappa < 0$, this asymptotic stability is a small-time stability due to Corollary \ref{cor2}. This concludes the proof of Theorem~\ref{thmcascaded}.
\hfill $\Box$

%\begin{remark}
%\label{remchain}
%To deal with chain of integrators, if the homogeneous stabilizing control law at a given step is not sufficiently smooth, one can approximate it by a differentiable homogeneous feedback law, to continue the application of Theorem \ref{thmcascaded} to stabilize the whole system.
%\end{remark}

%\Wilfrid{	
	%Let us introduce the $r$-homogeneous norm:
	%\begin{equation} \label{defi:hom_nom}
	%%\end{equation}

%	Moreover, the origin of $\dot x = f(x, \bar y (x))$ is supposed to be asymptotically stable and the vector-field $f$ to be $r$-homogeneous of degree $\kappa$, thus using \cite{Rosier:92:SCL} there exists $V(x)$ a $C^1$ , $\bar r$-homogeneous Lyapunov function of degree $d>-\kappa$ such that:
	
%\displaystyle \frac{\partial V}{\partial x} f(x, \bar y (x)) \leq - a V(x)^\frac{d+\kappa}{d}
%\label{eq:Vhombis},
%\end{equation}	
%where $a=-\sup_{\Vert x\Vert_{r}=1}  \frac{\partial V}{\partial x} f(x, \bar y (x))>0$.

%Since $f$ is $\bar r$-homogeneous, and $V$ is $\bar r$-homogeneous with respect to the weight $\bar r= (r_1, \cdots , r_{n} ) $   we conclude that $\frac{\partial V}{\partial x} f(x,y)$ is also $\bar r$-homogeneous and thus H\"older that is
%\[
%\left\Vert \frac{\partial V}{\partial x} f(x,y)-\frac{\partial V}{\partial x} f(x,\bar y (x)) \right\Vert \le b(x)
%\]

A natural question arise when $f$, the vector field involved in the cascade form (\ref{eqcascaded}), does not satisfy the is not homogeneous (i.e. does nor satisfy the homogeneity property \eqref{eq:prophom1}) but admits a homogeneous approximation $f_0$ at the origin. %OK pour Modif JM
Local homogeneous approximations has been investigated by Massera \cite{Massera56}, Rosier \cite{1992-Rosier-SCL}, Andrieu \cite{Andrieu08}, and recently Efimov \cite{Efimov15}. Let $g: \R^k \rightarrow \R^k$ be a continuous vector field. Let $ r=(r_1,\ldots,r_k)\tr \in (0,+\infty)^k$. Let us assume that $g$ has a $r$-homogeneous approximation of  at the origin, i.e. that there exists $g_0 : \R^k \rightarrow \R^k$	such that   for any $x \in S_{r}=\{x\in \R^k:  \sum_{i=1}^n|x_i|^{1/r_i}=1\}$ we  have $\lim_{\lambda \rightarrow 0^+} \lambda^{-\kappa}\Lambda_{r} (\lambda^{-1},g(\Lambda_{r} (\lambda, x)))=g_0(x)$  uniformly  on $S_{\bar r}$ (for some $\kappa \ge - \min_{1\le i\le k} r_i$, $\kappa$ is called the homogeneous degree of $g_0$). One has the following theorem
\begin{theorem}
\label{th-robustness-stability}
 If $0\in \R^k$ is locally asymptotically stable for
 $\dot  x= g_0(x)$ then $0\in \R^k$ is locally asymptotically stable for $\dot x = g(x)$. Moreover
 if $\kappa <0$ and if $0\in \R^k$ is small-time stable for $\dot  x= g_0(x)$ then $0\in \R^k$ is small-time stable for $\dot x = g(x)$
\end{theorem}
The  asymptotical statement of this theorem is \cite[Theorem 3]{1992-Rosier-SCL}. The finite-stability statement follows from  the proof of  \cite[Theorem 3]{1992-Rosier-SCL} given in \cite{1992-Rosier-SCL}. As a consequence of this robustness of the small-time stability we have the following corollary.
\begin{corollary}\label{corcascaded}
 Let $f_0:\R^n\times\R \to  \R^n$  which satisfies all hypothesis of theorem \ref{thmcascaded}: $\kappa\in \R$ satisfies (\ref{kappagd}), there exist $l\in (0,+\infty)$ and a continuous feedback law
		$\bar y : \mathbb{R}^n \rightarrow \mathbb{R}$  such that (\ref{lgrand}- \ref{homog-bary}- \ref{asforbary}- \ref{new-assumption-l}) holds. Let $f:\R^n\times\R \to  \R^n$ be a continuous map such that, uniformly on
$\Sigma:=\{(x,y)\in \R^n\times \R; \, |y|^{1/r_{n+1}}+\sum_{i=1}^n|x_i|^{1/r_i}=1\}$,
\begin{equation}
\label{approx-hom}
\lim_{\lambda \rightarrow 0^+} \lambda^{-\kappa}\Lambda_{\bar r}(\lambda^{-1}, f(\Lambda_{\bar r} (\lambda, x), \lambda^{r_{n+1}} y ))=f_0(x,y).
\end{equation}
Then the control law given by (\ref{feedback-v(x,y)}), with a sufficiently high gain $k>0$, makes $0\in \R^{n+1}$  locally asymptotically stable for the cascaded system (\ref{eqcascaded},$f$). Moreover, if $\kappa < 0$, this local asymptotic stability is a small-time stability.
\end{corollary}\noindent
For the proof of this corollary, it suffices to apply the previously obtained result of \cite{1992-Rosier-SCL} to
$k:=n+1$, $\R^{n+1}\simeq \R^n\times \R$, $g(x,y):=(f(x,y), v(x,y))$, and $g_0(x,y):=(f_0(x,y), v(x,y))$.  In fact, since our approach also gives
 homogeneous Lyapunov functions, one can provides a direct proof of Corollary~\ref{corcascaded}. Indeed,  using $W(x,y)$ given by (\ref{eq:Whom}), we have $
%\begin{equation}\label{eq:derWcorollary}
\left. \dot W\right |_{(\ref{eqcascaded},f)}=\left. \dot W\right |_{(\ref{eqcascaded},f_0)}+\frac{\partial W}{\partial x}(f-f_0).$
%\end{equation}
Since $W$ is $r$-homogeneous of degree $\alpha$ and its time derivative along the vector field $f_0$ with feedback (\ref{feedback-v(x,y)}) is $r$-homogeneous of degree $\alpha+\kappa$,  we have (see \cite{Moulay06})
$\left. \dot W\right |_{(\ref{eqcascaded},f_0)} \le - b_0 W^{\frac{\alpha+\kappa}{\alpha}},$ (see (\ref{Wpointfin}))
for some positive constant $b_0$ (take $b_0=-\max_{(x,y)\in \Sigma}\left(\frac{\partial W}{\partial x}f_0(x,y)+\frac{\partial W}{\partial y}v\right) 
% OK \modifsjm{W^{-(\alpha + \kappa)/\alpha}}
W^{-(\alpha + \kappa)/\alpha}$). Since $f_0$ is the homogeneous approximation at 0 of $f$ in the sense that \eqref{approx-hom} holds and $\frac{\partial W}{\partial x_i}$ is $r$-homogeneous (of degree $\alpha-r_i$) it implies that $\frac{\partial W}{\partial x}(f-f_0)(\Lambda_{\bar r}(\lambda, x), \lambda^{r_{n+1}}y)
% OK \modifsjm{\lambda^{r_{n+1}}y}) 
=o(\lambda^{\alpha+\kappa})$ uniformly on $\Sigma$ as $\lambda\rightarrow 0^+$. Thus, for every $(x,y)\in \R^n\times \R$ with $|x|+|y|$ small enough,  we have:
\begin{equation}
\label{dotWleq}
\left \vert \frac{\partial W}{\partial x}(f-f_0)\right \vert \le \frac{b_0}{2} W^{\frac{\alpha+\kappa}{\alpha}}.
\end{equation}
Hence,  for every $(x,y)\in \R^n\times \R$ with $|x|+|y|$ small enough,
\begin{equation}
\left. \dot W\right |_{(\ref{eqcascaded},f)}\le - \frac{b_0}{2} W^\frac{\alpha+\kappa}{\alpha},
\end{equation}
which whows that $0\in \R^{n+1}$  locally asymptotically stable for the cascaded system (\ref{eqcascaded},$f$). When $\kappa <0$, we have $(\alpha+\kappa)/\alpha<1$  and thus concludes from \eqref{dotWleq} that $0\in \R^{n+1}$ is small-time stable for this cascaded system (\ref{eqcascaded},$f$).

%\begin{lemma}\label{lem:ineq3} Let $a_0,a_1,a_2>0$ and $\beta\ge1-\kappa,\kappa<0$, then $\forall z\in\R_+ : \psi(z)=-a_0-a_1z^{\beta}+a_2z^{\beta+\kappa}$ is strictly negative iff $a_2^{-\frac{\beta}{\kappa}}(1+\frac{\kappa}{\beta})^{-(1+\frac{\beta}{\kappa})})\le -\frac{\beta}{\kappa}a_0a_1^{-(1+\frac{\beta}{\kappa})}$.
%\end{lemma}
%\begin{proof}
%	$\psi^\prime(z)=-a_1\beta z^{\beta-1}+a_2(\beta+\kappa)z^{\beta-1+\kappa}$ which zeros at $0$ and $z_0^{\kappa}=\frac{a_1\beta}{a_2(\beta+\kappa)}$ thus $\psi$ on $[0,+\infty[$ reach a maximum at $z_0$. We conclude noticing that $\psi(z_0)=-a_0-a_1\left(\frac{a_1\beta}{a_2(\beta+\kappa)}\right)^{\frac{\beta}{\kappa}}+a_2\left(\frac{a_1\beta}{a_2(\beta+\kappa)}\right)^{1+\frac{\beta}{\kappa}}$  is negative when $a_2^{-\frac{\beta}{\kappa}}(1+\frac{\kappa}{\beta})^{-(1+\frac{\beta}{\kappa})})\le -\frac{\beta}{\kappa}a_0a_1^{-(1+\frac{\beta}{\kappa})}$.
%\end{proof}
%}

\subsection{Small-time stabilization of the double integrator}
%The origin of the double integrator system (\ref{DI}) is GFTS with the following feedback law:
%\begin{equation}\label{FLDI}
%u(x_1,x_2)=-k_{2}\{\{ x_{2}\}^{\frac{1}{1+\kappa}}+k_1x_1\}^{1+2\kappa},
%\end{equation}
%where $-\frac{1}{2}<\kappa<0$; $k_1,k_2$ are strictly positive constants and
%$\{x\}^{\alpha}={sgn}(x)|x|^{\alpha}$, for any real number $\alpha\geq0$ and for all $x\in\mathbb{R}$.
%\end{theorem}
%The proof is based on the following Lyapunov function:
%\begin{equation}
%	\label{LyaDI}
%V(x_1,x_2)= \frac{1+k_1^{\frac{2+\kappa}{1+\kappa}}}{2+\kappa} |x_1|^{2+\kappa}+ \frac{1+\kappa}{2+\kappa}|x_2|^{\frac{2+\kappa}{1+\kappa}}+ k_1^{\frac{1}{1+\kappa}} x_1x_2
%\end{equation}

%
%\begin{theorem}\label{thmDI} The origin of the double integrator system (\ref{DI}) is GFTS with the following feedback law:
%\begin{equation}
%\label{FLDI}
%u = - k_1 sgn(x_1) |x_1 \mid^{\frac{\alpha}{2-\alpha}} - k_2 sgn(x_2) |x_2 \mid^{\alpha}
%\end{equation}
%where $k_1$ and $k_2$ are strictly positive constants.
%\end{theorem}
%The proof is based on the following Lyapunov function:
%\begin{equation}
%\label{LyaDI}
%V (x) = \frac{k_1 (2- \alpha)}{2} sgn(x_1) |x_1 \mid^{\frac{2}{2-\alpha}}+ \frac{x_2 ^2}{2}
%\end{equation}
%the conclusion holds applying LaSalle invariance principle showing that the closed-loop system is GAS. Moreover, the closed-loop system being homogeneous of negative degree, the GFTS of the origin is obtained using a result by Rosier \cite{RosierThese} (see \cite{bernuau-perruquetti} for details).
Let us consider the double integrator system:
\begin{equation}
\label{DI}
\dot x_1 = x_2, \;
\dot x_2 = u,
\end{equation}
where the state is $x = (x_1 , x_2 )^{\top}\in \mathbb{R}^2$ and the control is $u\in \mathbb{R}$.

Some results exist concerning the small-time stabilization of system (\ref{DI}). One can cite for example \cite{1986-Haimo-SICON, 1991-Coron-Praly-SCL, 2000-Bhat-Bernstein-SICON} for bounded continuous time-invariant small-time stabilizing feedback laws or \cite{bernuau-perruquetti} for an output-feedback control design. These feedback laws are singular when $x_1=0$ or $x_2=0$, therefore in order to apply ``desingularizing functions'' as explained in \cite{1991-Praly-dAndrea-Novel-Coron-IEEE} for cascaded system, and similarly in Theorem \ref{thmcascaded}, we refer to the following result inspired by \cite{2002-Hong-SCL} and \cite{bernuau-perruquetti}.

\begin{theorem}\label{thmG1} Let $-\frac{1}{2}<\kappa<0$ and $k_1>0,k_2>0$ be given. Then $0\in \R^2$ is small-time stable for the double integrator system (\ref{DI}) if one uses the following feedback law:
	\begin{equation}
		\label{FLDIBP}
		u(x_1,x_2) = - k_1 \{ x_1 \}^{\left\{1+2\kappa\right\}} - k_2 \{ x_2 \}^{\left\{\frac{1+2\kappa}{1+\kappa}\right\}}.
	\end{equation}
\end{theorem}
\begin{remark} The feedback law \eqref{FLDIBP} is a limiting case for \cite{1986-Haimo-SICON} which is not covered by \cite{1986-Haimo-SICON}.
\end{remark}
\noindent
{\bf Proof.} The time derivative along the solutions of the closed loop system (\ref{DI})-(\ref{FLDIBP}) of the following Lyapunov function:
%\begin{equation}
%\label{LyaNew}
$V (x_1,x_2) = \frac{k_1}{2(1+\kappa)} |x_1|^{2(1+\kappa)}+ \displaystyle \frac{x_2 ^2}{2}$,
%\end{equation}
 is $\dot V=-k_2 |x_2|^{\frac{2+3\kappa}{1+\kappa}} \le 0$. The LaSalle
 invariance principle shows that the closed-loop system is globally asymptotically stable. Moreover, the closed-loop system
 being % ok \modifsjm{$(1,1+\kappa)$-}
 $(1,1+\kappa)$-homogeneous of negative degree (this degree is $\kappa$), the small-time stability
 of the origin follows from Corollary~\ref{cor2}. % OK pour \modifsjm{
 In fact, instead of using the LaSalle invariance principle, one can construct an homogeneous strict Lyapunov, which is interesting for robustness issues. Indeed, for $\varepsilon>0$, let us consider
\begin{equation}
 \label{defVepsilon}
V_{\varepsilon}(x_1,x_2):= \frac{k_1}{2(1+\kappa)} |x_1|^{2(1+\kappa)}+ \displaystyle \frac{x_2 ^2}{2}
 +\varepsilon x_2\pui{x_1}{1+\kappa},
\end{equation}
which is $(1,1+\kappa)$-homogeneous of degree $2(1+\kappa)$. Then straightforward computations give the existence of an $\varepsilon_0>0$ such that,
for every $\varepsilon\in (0,\varepsilon_0)$, $(1/2)V\leq V_{\varepsilon}\leq 2 V$ and there exists $c>0$ (depending on $\varepsilon$) such that $\dot V_{\varepsilon}\leq -c V_{\varepsilon}^{\frac{2+3\kappa}{2(1+\kappa)}}$
\hfill $\Box$

% OK enlevé \commentsjm{je propose d'enlever le th\'{e}or\`{e}me qui suit (Theorem \ref{thmG2}) car il est un cas particulier du Theorem \ref{thmDI} : prendre $l=(1-\kappa)/(1+\kappa)$. De plus consid\'{e}rer $l=1/(1+\kappa)$ au lieu de  $l=(1-\kappa)/(1+\kappa)$ aurait \'{e}t\'{e} plus naturel et donne effectivement des formules plus simples. C'est ce $l$ que j'ai choisi pour le slider plus loin.
\begin{theorem}\label{thmG2} Let
\begin{equation}
-\frac{1}{2}<\kappa<0,\, k_1>0,\, k_2> \frac{1-\kappa^2}{2+\kappa} 2^{\frac{1-4\kappa -3\kappa^2}{1-\kappa^2}}k_1^{\frac{1}{1+\kappa}}
\end{equation}
be given. % OK vérifié \commentsjm{avant, pour $k_2$ il y avait $k_2> 2^{\frac{1+\kappa^2}{1-\kappa^2}}k_1^{\frac{1}{1+\kappa}}$.} 
Then $0\in \R^2$ is small-time stable for the double integrator system (\ref{DI}) if one uses the following feedback law:
	\begin{equation}\label{FLDIH}
	u(x_1,x_2)=-k_{2}\pui{ \pui{x_2}{\frac{1-\kappa}{1+\kappa}}+ k_1^{\frac{1-\kappa}{1+\kappa}} \pui{x_1}{1-\kappa} }{\frac{1+2\kappa}{1-\kappa}}.
		%k_{2}\left\{\{ x_{2}\}^{\left\{\frac{1-\kappa}{1+\kappa}\right\}}+\{k_1 \{ x_1\}^{\left\{1+\kappa\right\}}\}\right\}^{\{\frac{1+2\kappa}{1-\kappa}\}}
	\end{equation}
\end{theorem}
\noindent
{\bf Proof (sketch).} Let us introduce %the following functions
	$\bar x_{2}=-k_1 \{ x_1\}^{\left\{1+\kappa\right\}},%\phi_{1}=\left|x_1\right|,
	\varphi_{2}(x_1,y)=\{y\}^{\left\{\frac{1-\kappa}{1+\kappa}\right\}}-\{\bar x_{2}\}^{\left\{\frac{1-\kappa}{1+\kappa}\right\}}
	$  and
	$\Phi_{2}(x_1,x_2)= \int_{\bar x_2}^{x_2} \varphi_{2}(x_1,y)dy
	= \frac{1+\kappa}{2}\left(\left|x_{2}\right|^\frac{2}{1+\kappa}+\frac{1-\kappa}{1+\kappa}\left|\bar x_{2}\right|^\frac{2}{1+\kappa}\right)-x_{2}\{\bar x_{2}\}^{\frac{1-\kappa}{1+\kappa}}$
	which is nonnegative and zero if and only if $x_2=\bar x_{2}(x_1)$. Thus we have $	u(x_1,x_2)=-k_{2} \pui{	\varphi_{2}(x_1,x_2)}{\frac{1+2\kappa}{1-\kappa}}$.
The time derivative along the solutions of the closed loop system (\ref{DI}) (\ref{FLDIH}) of the following Lyapunov function:
\begin{equation}
	\label{LyaDIH}
	W(x_1,x_2)= \frac{x_1^{2}}{2}+a
	 \Phi_{2}(x_1,x_2),
\end{equation}
is $\dot W=x_1(x_2+\bar x_{2}-\bar x_{2})+a \dot \Phi_{2}=Q_1+Q_2+a(Q_3+Q_4)$, where $Q_1=x_1\bar x_{2}= -k_1\vert x_1\vert^{2+\kappa}, Q_2=x_1(x_2-\bar x_{2}) ,
Q_3=\frac{\partial \Phi_{2}}{\partial x_1}x_2= (1-\kappa) k_1^{\frac{1-\kappa}{1+\kappa}} \vert x_1\vert^{-\kappa}(x_2-\bar x_{2})x_2,
Q_4=\frac{\partial \Phi_{2}}{\partial x_2}u=\varphi_{2}(x_1,x_2)u(x_1,x_2)=-k_2 \vert \varphi_{2} (x_1,x_2)\vert^{\frac{2+\kappa}{1-\kappa}}$. Selecting $a(1-\kappa)
k_1^{\frac{2}{1+\kappa}}=1$ we get
$\dot W=-k_1\vert x_1\vert^{2+\kappa}-ak_2 \vert \varphi_{2}(x_1,x_2) \vert^{\frac{2+\kappa}{1-\kappa}}+\frac{1}{k_1} \vert x_1\vert^{-\kappa} (\bar x_{2}-x_2)^2$. Using Lemma \ref{lem:ineq1} with $l:=(1-\kappa)/(1+\kappa)$ we get $(\bar x_{2}-x_2)^2
\le 2^{\frac{-4\kappa}{1-\kappa}}  \vert \varphi_{2}(x_1,x_2)
 \vert^{\frac{2(1+\kappa)}{1-\kappa}}$  and thus $\dot W \le -k_1\vert x_1\vert^{2+\kappa}-ak_2 \vert \varphi_{2}(x_1,x_2) \vert^{\frac{2+\kappa}{1-\kappa}}+\frac{2^{\frac{-4\kappa}{1-\kappa}}}{k_1}
\vert x_1\vert^{-\kappa}  \vert \varphi_{2}(x_1,x_2) \vert^{\frac{2(1+\kappa)}{1-\kappa}}$. Clearly $\dot W$ is negative when $\vert x_1\vert=0$. For $x_1\neq0$ the signum of the right-hand side of this inequality is the same as the signum of $\psi(z)=-k_1-ak_2 z^{2+\kappa}+\frac{2^{\frac{-4\kappa}{1-\kappa}}}{k_1}z^{2(1+\kappa)}$
(consider $z:=\vert \varphi_{2}(x_1,x_2)\vert ^{1/(1-k)}/|x_1|$).  Lemma \ref{lem:ineq2} (with $\beta=2+\kappa>\gamma=2(1+\kappa)\ge 1,a_0=k_1,a_1=ak_2,a_2=\frac{2^{\frac{-4\kappa}{1-\kappa}}}{k_1}$) and
  and Remark~\ref{rem.simplecondition} show that $\dot W$ is negative definite under the given conditions for $\kappa, k_1,k_2$. The closed-loop system being homogeneous of negative degree, the small-time stability of the origin follows from Corollary~\ref{cor2}.
\hfill $\Box$

% OK pour \modifsjm{We can also apply our Theorem~\ref{thmcascaded} to stabilize in small time $0\in \R^2$ for system (\ref{DI}). As already mentioned in Remark~\ref{rem-Morin-Samson} it gives new feedback laws which are more suitable than the ones given in Theorem~\ref{thmG1} when one wants to add ``integrators'', as we are going to do in Section~\ref{sec-slider} in order to treat the slider control system. Moreover using our proof of Theorem \ref{thmcascaded} we can get precise lower bounds on the gain for the feedback laws as shown in the following theorem.}
We can also apply our Theorem~\ref{thmcascaded} to stabilize in small time $0\in \R^2$ for system (\ref{DI}). As already mentioned in Remark~\ref{rem-Morin-Samson} it gives new feedback laws which are more suitable than the ones given in Theorem~\ref{thmG1} when one wants to add ``integrators'', as we are going to do in Section~\ref{sec-slider} in order to treat the slider control system. Moreover using our proof of Theorem \ref{thmcascaded} we can get precise lower bounds on the gain for the feedback laws as shown in the following theorem.

\begin{theorem}\label{thmDI}
Let $-1/2<\kappa<0$, $l \ge 1/(1+\kappa)$,
% \modifsjm{$l \ge 1/(1+\kappa)$,}  OUI ca marche \commentsjm{avant on avait $l\geq 2$ ; v\'{e}rifier que la nouvelle condition, plus faible, est bien suffisante.} 
and $k_1>0, k_2>g(l,\kappa,k_1)$ be given, where
\begin{equation}\label{eq:condition_gain}
g(l,\kappa,k_1):=
% OK VÉRIFIÉ \modifsjm{\frac{2l(1+\kappa)^2}{(l+1)(1+\kappa)+\kappa}}
\frac{2l(1+\kappa)^2}{(l+1)(1+\kappa)+\kappa}
2^{\frac{(l-1)\left((l+1)(1+\kappa)+\kappa\right)}{l(1+\kappa)}}
k_1^{\frac{1}{1+\kappa}}.
\end{equation}
% OK vérifié \commentsjm{avant le terme en rouge \'{e}tait $l(1+\kappa)$ et il y avait $\left(\frac{(l+1)(1+\kappa)-\kappa-2}{(l+1)(1+\kappa)+\kappa}\right)^{\frac{(l+1)(1+\kappa)-\kappa-2}{2(1+\kappa)}}$} 
Then %there exists a sufficiently high gain $k_2>0$ such that
$0\in \R^2$ is small-time stable for the double integrator system (\ref{DI}) if one uses the following feedback law:
\begin{equation}
\label{FLDI}
u(x_1,x_2) = - k_2 \pui{\pui{x_2}{l} - \pui{\bar x _2}{l}}{\frac{1+2\kappa}{l(1+\kappa)}},
\end{equation}
where
\begin{equation}
\label{FLDI1}
\bar x_2 = -k_1 \pui{x_1}{1+\kappa}.
\end{equation}
\end{theorem}

\noindent
{\bf Proof.} Let us introduce $(r_1,r_2)=(1, 1+ \kappa)$ as the weight of the dilation. The proof is a straightforward application of Theorem \ref{thmcascaded}, noticing that $0\in \R$ is small-time stable for
\begin{equation}
\dot x_1 = - k_1 \pui{x_1}{1+\kappa}=:\bar x_2 (x_1),
\end{equation}
the control law $\bar x_2 (x_1)$ being $r_1$-homogeneous of degree $r_2$, since $r_2 = (1+\kappa)r_1$.  It is clear that, since
$l\geq 1/(1+\kappa)$,  $\left\{\bar x _2\right\}^{\{l\}}(x_1)$ is $C^1$ and \eqref{lgrand} holds. Moreover, the vector-field associated to subsystem $x_1$ is $(r_1, r_2)$-homogeneous of degree $\kappa$ (i.e. we have \eqref{eq:prophom1}). %, since from (\ref{FLDI2}) we have for any positive real number $\lambda$:
%$$
%\lambda ^{r_2} = \lambda^{\kappa_1 + r_1} .
%$$
The closed-loop system (\ref{DI}) with the feedback law (\ref{FLDI}) is itself $(r_1, r_2)$-homogeneous of degree $\kappa$. We have $\kappa + r_2 > 0$ (because $r_2=1+\kappa$ and $-\frac{1}{2}<\kappa$ ), so that condition (\ref{kappagd}) of Theorem \ref{thmcascaded} holds, meaning that the control law (\ref{FLDI}) is well defined and continuous. Since $(l+1)>1+\kappa$ we also have \eqref{lgrand}. Lastly, $\kappa<0$ which implies that $0\in \R^2$ is small-time stable for the double integrator system (\ref{DI})  if $k_2>0$ is large enough. In fact, following the proof of Theorem \ref{thmcascaded} one can get the lower bound for $k_2$ given in Theorem~\ref{thmDI}.  Let $\alpha:=(l+1)(1+\kappa)>1$. Setting $\varphi(x_1,s)=\pui{s}{l} - \pui{\bar x_2}{l}$ we have $u(x_1,x_2)=-k_{2} \pui{	\varphi(x_1,x_2)}{\frac{1+2\kappa}{l(1+\kappa)}}$. Let  $W(x_1,x_2)=\frac{\vert x_1 \vert^{\alpha}}{\alpha}+a\Phi (x_1,x_2),$ where $\Phi (x_1,x_2)=\int _{\bar x_2} ^{x_2} \varphi(x_1,s)ds$. %=\int _{\bar x_2} ^{x_2} (\pui{s}{l} - \pui{\bar x_2}{l}) ds
%	= \displaystyle \frac{l}{l+1} \vert \bar x_2\vert^{l+1} (x_1) - x_2 \pui{\bar x_2}{l} + \displaystyle \frac{\vert x_2\vert ^{l+1}}{l+1}$ (which is positive and zero only when $x_2=\bar x_{2}(x_1)$).
Its time derivative along the solution of the closed-loop system
is $\dot W=\pui{x_1}{\alpha-1}(\bar x_{2}+x_2-\bar x_{2})+a \dot \Phi=Q_1+Q_2+a(Q_3+Q_4)$, where $Q_1=\pui{x_1}{\alpha-1}\bar x_{2}= -k_1\vert x_1\vert^{\alpha+\kappa}, Q_2=\pui{x_1}{\alpha-1}(x_2-\bar x_{2}) , Q_3=\frac{\partial \Phi}{\partial x_1}x_2= l (1+\kappa)k_1^{l} \vert x_1\vert^{\alpha-\kappa-2}(x_2-\bar x_{2})x_2,Q_4=\frac{\partial \Phi}{\partial x_2}u=-k_2
\vert \varphi \vert^{\frac{\alpha+\kappa}{l(1+\kappa)}}$.
Selecting $al(1+\kappa)k_1^{1+l}=1$ we get that $\dot W=-k_1\vert x_1\vert^{\alpha+\kappa}-ak_2 \vert \varphi (x_1,x_2)  \vert^{\frac{\alpha+\kappa}{l(1+\kappa)}}+\frac{1}{k_1}\vert x_1\vert^{\alpha-\kappa-2} (\bar x_{2}-x_2)^2$. Using Lemma \ref{lem:ineq1} we get $(\bar x_{2}-x_2)^2 \le 2^{2\frac{l-1}{l}} \vert \varphi(x_1,x_2) \vert^{\frac{2}{l}}$ and  thus $\dot W \le -k_1\vert x_1\vert^{\alpha+\kappa}-ak_2 \vert \varphi(x_1,x_2) \vert^{\frac{\alpha+\kappa}{l(1+\kappa)}}+\frac{2^{2\frac{l-1}{l}}}{k_1} \vert x_1\vert^{\alpha-\kappa-2} \vert \varphi(x_1,x_2) \vert^{\frac{2}{l}}$. Clearly $\dot W$ is negative when $\vert x_1\vert=0$. For $x_1\neq0$ the signum of the right-hand side of this inequality is the same as the signum of $\psi(z)=-k_1-ak_2 z^{\alpha+\kappa}+\frac{2^{2\frac{l-1}{l}}}{k_1} z^{2(1+\kappa)}$.
Lemma \ref{lem:ineq2} and Remark~\ref{rem.simplecondition}(with $\beta=\alpha+\kappa>\gamma=2(1+\kappa)>1,a_0=k_1,a_1=ak_2,a_2=\frac{2^{2\frac{l-1}{l}}}{k_1} $) show that $\dot W$ is negative definite under the given conditions for $\kappa, k_1,k_2,l$. The closed-loop system being homogeneous of negative degree, the small-time stability of the origin follows from Corollary~\ref{cor2}.
%\begin{equation}
%\dot V \le \left(-\frac{k_1}{2} +a\right) \vert x_1\vert^{4+ 5\kappa} + (c_2+ac_4-ak_2 2^{-2}) \vert x_2 - \bar x_2 (x_1)\vert^{(4+5 \kappa)/(1+\kappa)}
%\end{equation}
%for some positive $c_2,c_4$. Selecting $a=\frac{k_1}{4}$ and $k_2$ sufficiently large ($k_2>\frac{c_2+ac_4}{a2^{-2}}$),
%we have
%%\begin{equation}
%\begin{gather}
%V(x_1,x_2)=\frac{\vert x_1 \vert^{4(1+\kappa)}}{4(1+\kappa)}+\frac{k_1}{4}\int _{\bar x_2} ^{x_2} (s^3 - \bar x_2 ^3 ) ds,
%\\
%\dot V \le -c_5 V^{\frac{1+5\kappa/4}{1+\kappa}}
%\end{gather}
%for some positive $c_5(k_2)$ (a continuous function of $k_2$ which goes to infinity with $k_2$) and conclude that the origin is small-time stable.
%}
\hfill $\Box$

This last result can be extended to following ``modified double integrator'' system:
\begin{equation}
\label{DI_modif_nu}
\left \{ \begin{array}{l}
\dot x_1 = \vert x_1\vert^{\nu} x_2\\
\dot x_2 = u
\end{array}
\right .
\end{equation}
where $0 \le \nu< 1$. We have the following theorem.
\begin{theorem}
\label{thm:DI_modif_nu} Let  $\nu \in [0,1)$, $-(1-\nu)/2<\kappa<0$, $l >1/(1+\kappa -\nu)$,
% OK ca marche \commentsjm{avant on supposait $l\geq 2$ ; bien v\'{e}rifier que la nouvelle condition plus faible est suffisante}, 
$k_1>0$, and  $k_2>h(l,\kappa,\nu,k_1)$ be given, where
{\scriptsize
\begin{gather}\label{def.h}
h(l,\kappa,\nu,k_1):=
% OUI vérifié \modifsjm{\frac{2l(1+\kappa-\nu)^2}{(l+1)(1+\kappa-\nu)+\kappa}2^{\frac{(l-1)\left((l+1)(1+\kappa-\nu)+\kappa\right)}{l(1+\kappa-\nu)}}}
\frac{2l(1+\kappa-\nu)^2}{(l+1)(1+\kappa-\nu)+\kappa}2^{\frac{(l-1)\left((l+1)(1+\kappa-\nu)+\kappa\right)}{l(1+\kappa-\nu)}} k_1^{\frac{1}{1+\kappa-\nu}}.
\end{gather}
}
% Tu as bien raison j'ai vérifié \commentsjm{avant la condition sur $k_2$ \'{e}tait $k_2> g(l,\kappa-\nu,k_1)$. Mais j'ai l'impression que si certains $\kappa$ peuvent \^{e}tre remplac\'{e}s par $\kappa-\nu$ ce n'est pas le cas de tous les $\kappa$. \`{A} v\'{e}rifier.} 
Then %there exists a sufficiently high gain $k_2>0$ such that
	$0\in \R^2$ is small-time stable for the system (\ref{DI_modif_nu}) if one uses the following feedback law:

	\begin{equation}
\label{FLDImodif}
u(x_1,x_2) = - k_2 \pui{\pui{x_2}{l} -\pui{ \bar x _2}{l}}{\frac{1+2\kappa-\nu}{l(1+\kappa-\nu)}},
\end{equation}

where
\begin{equation}
\label{eq:FLDImdif}
\bar x_2 = -k_1 \pui{x_1}{1+\kappa-\nu}.
\end{equation}
\end{theorem}

\noindent
{\bf Proof.} The case $\nu=0$ is proven in Theorem \ref{thmDI}. Let us introduce $(r_1,r_2) =(1, 1+ \kappa-\nu)\in (0,+\infty)\times(0,+\infty)$ as the weight of the dilation. Clearly the system (\ref{DI_modif_nu}) satisfies \eqref{eq:prophom1}.
%Indeed we have:
The proof is a straightforward application of Theorem \ref{thmcascaded}, noticing that the origin is small-time stable for the subsystem:
\begin{equation}
\dot x_1 = \vert x_1\vert^{\nu} \bar x_2 (x_1)=- k_1 \pui{x_1}{1+\kappa},
\end{equation}
the control law $\bar x_2 (x_1)$ being $r_1$-homogeneous of degree $r_2$, since $r_2 = 1+\kappa-\nu$.  Since
$l (1+\kappa -\nu)>1$, $\bar x_2^l (x_1)$ is of class $C^1$. Moreover, the vector-field associated to subsystem $x_1$ is $(r_1, r_2)$-homogeneous of degree $\kappa$. %, since from (\ref{FLDI2}) we have for any positive real number $\lambda$:
%$$
%\lambda ^{r_2} = \lambda^{\kappa_1 + r_1} .
%$$
Finally, the closed-loop system (\ref{DI}) with the feedback law (\ref{FLDImodif}) is itself $(r_1, r_2)$-homogeneous of degree $\kappa$. We have $\kappa + r_2 > 0$ (because $r_2=1+\kappa-\nu$ and $-(1-\nu)/2<\kappa$), so that condition (\ref{kappagd}) of Theorem \ref{thmcascaded} holds, meaning that the control law (\ref{FLDImodif}) is well defined and continuous. Lastly, $\kappa<0$ which implies that $0\in \R^2$ is small-time stable for the modified double integrator system given by (\ref{DI_modif_nu}) if $k_2>0$ is large enough.  In fact, as for Theorem~\ref{thmDI}, the proof of %Theorem \ref{thmcascaded} (see the proof of
Theorem~\ref{thmcascaded} allows to quantify ``$k_2>0$ large enough''.  We set $ \alpha=(l+1)(1+\kappa-\nu)>1$ and
%\begin{equation}
%V(x_1)=\frac{\vert x_1 \vert^{2q(1+\kappa-\nu)}}{2q(1+\kappa-\nu)}.
%\end{equation}
%Thus we
we consider the Lyapunov function
%\begin{equation}
$W(x_1,x_2)=\frac{\vert x_1 \vert^{\alpha}}{\alpha}+a\Phi (x_1,x_2)$,
%\end{equation}
where
%\begin{equation}
$\Phi (x_1,x_2)= \int _{\bar x_2} ^{x_2} (\pui{s}{l} - \pui{\bar x_2 }{l}) ds.$ % = \displaystyle \frac{2q-1}{2q} \bar x_2^{2q} (x_1) - x_2 \pui{\bar x_2}{2q-1}(x_1) + \displaystyle \frac{x_2^{2q}}{2q}.
%\end{equation}
Similarly to the proof of Theorem \ref{thmDI} we get $\dot W=\pui{x_1}{\alpha-1} \vert x_1\vert^\nu(\bar x_{2}+x_2-\bar x_{2})+a \dot \Phi=Q_1+Q_2+a(Q_3+Q_4)$, where $Q_1=\pui{x_1}{\alpha+\nu-1}\bar x_{2}= -k_1\vert x_1\vert^{\alpha+\kappa}, Q_2=\pui{x_1}{\alpha+\nu-1}(x_2-\bar x_{2}) , Q_3=\frac{\partial \Phi}{\partial x_1}x_2= l (1+\kappa-\nu)k_1^{l} \vert x_1\vert^{\alpha-\kappa-2+2\nu}(x_2-\bar x_{2})x_2,Q_4=\frac{\partial \Phi}{\partial x_2}u=-k_2 \vert \varphi \vert^{\frac{\alpha+\kappa}{l(1+\kappa-\nu)}}$, with $\varphi(x_1,x_2):=\pui{x_1}{l} - \pui{\bar x_2}{l}$.
Selecting $al(1+\kappa-\nu)k_1^{1+l}=1$ we get that $\dot W=-k_1\vert x_1\vert^{\alpha+\kappa}-ak_2 \vert \varphi (x_1,x_2)  \vert^{\frac{\alpha+\kappa}{l(1+\kappa-\nu)}}+\frac{1}{k_1}\vert x_1\vert^{\alpha-\kappa-2+2\nu} (\bar x_{2}-x_2)^2$. Using Lemma \ref{lem:ineq1} and then  Lemma \ref{lem:ineq2} (with $\beta=\alpha+\kappa>\gamma=2(1+\kappa-\nu)>1,a_0=k_1,a_1=ak_2,a_2=\frac{2^{2\frac{l-1}{l}}}{k_1} $) together with Remark~\ref{rem.simplecondition}  we conclude that $\dot W < 0$ in $\R^2\setminus\{0\}$ under the above given conditions for $l$, $k_1$, $k_2$, $\kappa$ and $\nu$. % time-derivative of $W$ along the solution of the closed-loop system satisfies
%\begin{equation}
%\dot V \le \left(-\frac{k_1}{2} +a\right) \vert x_1\vert^{2q(1+\kappa-\nu)+\kappa} + (c_2+ac_4-ak_2 2^{-2}) \vert x_2 - \bar x_2 (x_1)\vert^{\frac{2q(1+\kappa-\nu)+\kappa}{1+\kappa-\nu}}
%\end{equation}
%for some positive $c_2$ and $c_4$. Selecting $a=k_1/4$ and $k_2$ sufficiently large ($k_2>4(c_2+ac_4)/a$), we have
%\begin{equation}
%\dot V \le -c_5(k_2) V^{\frac{2q(1+\kappa-\nu)+\kappa}{2q(1+\kappa-\nu)}}
%\end{equation}
%for some positive $c_5(k_2)$ and
Finally using Corollary~\ref{cor2}, we conclude that $0\in \R^2$ is small-time stable.
%FTS with  a settling-time less than:
%\begin{equation}
%T \le \frac{2q(1+\kappa-\nu) V^{\frac{-\kappa}{2q(1+\kappa-\nu)}}(0)}{-\kappa c_5(k_2)},
%\end{equation}
%thus there exists a positive $c(k_2)$ such that the settling-time is less than
%\begin{equation}
%T \le c(k_2) \left(\vert x_1(0)\vert^{-\kappa}+\vert x_2(0)\vert^{\frac{-\kappa}{1+\kappa-\nu}} \right).
%\end{equation}
%Similar arguments as the ones used in the proof of Theorem \ref{thmDI}  show that $c(k_2)$ is a continuous function of $k_2$ which goes to 0 when $k_2$ goes to infinity.
\hfill $\Box$

%\section{Finite time stabilization of an elementary mobile robot}
%
%\begin{equation}\label{eq-mobile robot}
%\dot x_1=u_1\sin(u_2), \; \dot x_2= u_1\cos(u_2), \;
%\end{equation}

\section{Small-time stabilization of the unicycle robot}
\label{sec-unicycle}
Let us define $(x_1,x_2,x_3)\tr \in \R^6$ and  $(u_1,u_2)\tr\in\R^2$ by
\begin{equation}\label{defxu}
x_1=:x,\; x_2:=y, \; x_3:=\psi,\;
u_1:=\nu_1,\; u_2:=\Omega.
\end{equation}
Then the unicycle control system \eqref{eq:unicycle} becomes the control system

\begin{equation}
\label{unicyclexcossin}
\dot{x}_1  =   u_1\cos(x_3),\;
\dot{x}_2 =   u_1\sin(x_3),\;
\dot{x}_3  =   u_2
\end{equation}
where the state is $(x_1,x_2,x_3)\tr \in \R^3$ and
the control is $(u_1,u_2)\tr \in\R^2$.

One first observes that if $x_2$ and $x_3$ vanish at some $t_0$ and if, after time $t_0$, $u_2$ is equal to $0$, then $x_2$ and $x_3$ remain equal to $0$. The idea is then to use
a first stationary feedback law to steer $(x_2,x_3)\tr$ to $0\in \R^2$ and then take $u_2=0$ and use $u_1$ to steer $x_1\in \R$ to $0$.
\begin{remark}
This strategy is inspired from \cite{1995-Coron-SICON}. In \cite{1995-Coron-SICON} a preliminary stationary feedback law is used to steer the control system on (a neighborhood of) a special curve which can be sent to $0$ by a special feedback law. This strategy  is  already used in \cite{2016-Guilleron-PhDthesis}. The main novelty of our approach compared to
\cite{2016-Guilleron-PhDthesis} is to use this strategy in the framework of homogeneous approximation systems and homogeneous feedback laws. This will allow us to add integrators on the controls and therefore deal with the slider control system \eqref{eq:sousactionneinertiel}.
\end{remark}

We then consider the following quadratic homogeneous approximation of the control system \eqref{unicyclexcossin} (see \cite{2017-Coron-Rivas-SICON})
\begin{equation}
\dot{x}_1  =  u_1,\;
\dot{x}_2 =  u_1 x_3,\;
\dot{x}_3  =   u_2,
\label{unicycle-simp1}
\end{equation}
which is a control system where the state is $(x_1,x_2,x_3)\tr \in \R^3$ and
the control is $(u_1,u_2)\tr \in\R^2$.

As for \eqref{unicyclexcossin}, let us observe that, for \eqref{unicycle-simp1}, if $x_2$ and $x_3$ vanish at some $t_0$ and if, after time $t_0$, $u_2$ is equal to $0$, then $x_2$ and $x_3$ remain equal to $0$. The idea is then, again, to use, for \eqref{unicycle-simp1},
a first stationary feedback law to steer $(x_2,x_3)\tr$ to $0\in \R^2$ and then take $u_2=0$ and use $u_1$ to steer $x_1\in \R$ to $0$.

 For this homogeneous approximation, a time-piecewise continuous periodic feedback law will be designed leading to the small-time stabilization of \eqref{unicycle-simp1} (Theorem \ref{thm:unicycle_LFTS_1}).  Then, using homogeneity arguments, we will check that the same feedback laws leads
to the small-time stabilization of the orignal model system \eqref{sliderxcossin}.

In this first design, the control is derived using Theorem \ref{thm:DI_modif_nu}. For this, in a first step, let us choose $u_1=\vert x_2\vert^\nu, \nu \in (0,1)$ and $\bar x_3=-k_2  \pui{x_2}{1+\kappa_2-\nu},  \kappa_2 \in (-(1-\nu)/2,0)$, so that the "equivalent" dynamics of the subsystem $x_2$ is
$ \dot x_2=-k_2  \pui{x_2}{1+\kappa_2}$ ($x_2$ converges to zero in small time) and then, the control $u_2$ is designed using Theorem \ref{thm:DI_modif_nu} (a consequence of our  cascading result Theorem \ref{thmcascaded}) : thus we obtain $x_2 (t) = x_3 (t) = 0$ for all $t \geq T/2$ provided that $x_2 (0)$ and $x_3 (0)$ are small enough. Note than on this time interval $x_1$ increase until $x_2$ reach zero because then $u_1=0$: we should just show that $x_1$ remains bounded on a finite time interval. Then, in a second step, note that setting $u_2 = 0$ for $t \in (T/2, T)$ then $x_3,x_2$ will remain zero (if they have reached zero). Thus in order to stabilize $x_1$ we can choose for $t \in (\frac{T}{2}, T)$, $u_1=-k_1 \pui{x_1 }{1+\kappa_1} , -1<\kappa_1<0$.
This leads to the following theorem:

\begin{theorem}\label{thm:unicycle_LFTS_1}
Let $T>0$, $\nu\in(0,1)$, $-1<\kappa_1<0$, $-(1-\nu)/2<\kappa_2<0$,  $l>1/(1+\kappa_2-\nu)$, $k_1>0$, and $k_2>0$. Then, for every $k_3>h(l,\kappa_2,\nu,k_2)$ %\Wilfrid{Avant on avait $\kappa$}
% OK je suis d'accord \modifsjm{for every $k_3>h(l,\kappa,\nu,k_2)$} \commentsjm{avant il y avait ``there exists large enough $k_3>h(l,\kappa,\nu,k_2)$'', mais cette formulation laisse sugg\'{e}rer qu'il est possible qu'il existe des $k_3>h(l,\kappa,\nu,k_2)$ qui ne conviennent pas. Ce n'est pas le cas n'est-ce pas ?}
(where the function $h$ is given by (\ref{def.h})) for the feedback law $u=(u_1,u_2)\in C_{\text{t-p}} (\mathbb{R} \times \mathbb{R}^3; \mathbb{R}^2)$ defined by
\begin{gather}
		u_1 = \vert x_2\vert^{\nu} \text{ in } \textstyle \left[0,\frac{T}{2} \right]\times \mathbb{R}^3,	\label{u1_FTS_unicycle_1}\\
\label{unicycle-u-2-0-T/2}
		u_2 = - k_3\pui{\pui{x_3}{l}-\pui{\bar x_3}{l}}{\frac{1+2\kappa_2-\nu }{l(1+\kappa_2-\nu)}}
			% OK \modifsjm{\pui{x_3}{l}} - \modifsjm{\pui{\bar x_3}{l}}}{\frac{1+2\kappa_2-\nu }{\modifsjm{l}(1+\kappa_2-\nu)}}
\text{ in } \textstyle \left[0,\frac{T}{2}\right]\times \mathbb{R}^3, %\\	
\text{ with }
		\bar x_3=-k_2  \pui{x_2}{1+\kappa_2-\nu},	\nonumber\\
		u_1 = - k_1 \pui{ x_1 }{1+\kappa_1},  \text{ in } \textstyle \left(\frac{T}{2},T\right)\times \mathbb{R}^3,\\
		u_2 = 0 \text{ in } \textstyle \left(\frac{T}{2},T\right)\times \mathbb{R}^3,\\
		u(t+T,x)=u(t,x),\; \forall t \in \mathbb{R},\; \forall x \in \mathbb{R}^3, 	\label{u1u2_unicycle_en0_fast4}
\end{gather}
 % OK je suis d'accord \commentsjm{avant il y avait $u_2=- k_3\pui{x_3 - \bar x_3}{\frac{1+2\kappa_2-\nu }{1+\kappa_2-\nu}}$ dans \eqref{unicycle-u-2-0-T/2}; c'est une possibilit\'{e}, mais alors on n'a pas \'{e}tabli de borne inf\'{e}rieure sur $k_3$ et cela ne rentre pas bien dans ce que l'on a fait avant}
 and for $f(x,u)= (u_1,\; u_1 x_3,\; u_2)^{\top}$, \eqref{unulen0}, \eqref{uperiodic}, and \eqref{0toujours} hold  and there exists $\epsilon >0$ such that \eqref{epsilontoujours} holds too.
 In particular  (\ref{unicycle-simp1}) is small-time locally stabilizable by means of explicit time-varying piecewise
 continuous % OK \modifsjm{stationary} 
 stationary feedback laws.
\end{theorem}

\noindent
{\bf Proof.}
On $\left[0,T/2\right]$, the closed-loop system is
\begin{equation}
\left\{
\begin{array}{lcl}
\dot{x}_1 & = &   \vert x_2\vert^{\nu},\\
\dot{x}_2 & = &   \vert x_2\vert^{\nu} x_3,\\
% OK pour les modifs de JM 
\dot{x}_3 & = & - k_3\pui{\pui{x_3}{l}-\pui{\bar x_3}{l}}{\frac{1+2\kappa_2-\nu }{l(1+\kappa_2-\nu)}}.
\end{array}
\right.
\label{cl1_unicycle-simp1}
\end{equation}
% OK pour les modifs de JM
We apply the proof of Theorem \ref{thm:DI_modif_nu} to the $(x_2,x_3)-$subsystem. We set $ \alpha=(l+1)(1+\kappa_2-\nu)$ and
%\begin{equation}
%V(x_1)=\frac{\vert x_1 \vert^{2q(1+\kappa-\nu)}}{2q(1+\kappa-\nu)}.
%\end{equation}
%Thus we
we consider the Lyapunov function
%\begin{equation}
$W(x_2,x_3)=\frac{\vert x_2 \vert^{\alpha}}{\alpha}+a\Phi (x_2,x_3)$,
%\end{equation}
where
%\begin{equation}
$\Phi (x_2,x_3)= \int _{\bar x_3} ^{x_3} (\pui{s}{l} - \pui{\bar x_3 }{l}) ds$ and $a$ is such that
$al(1+\kappa_2-\nu)k_2^{1+l}=1$. From the proof of Theorem \ref{thm:DI_modif_nu} we know that
\begin{equation}\label{decayW}
\dot W\leq -\frac{1}{C}W^{\frac{\alpha+\kappa_2}{\alpha}}. 
\end{equation}
%\Wilfrid{Avant on avait $W^{\frac{\alpha_2+\kappa_2}{\alpha_2}}$ }
In \eqref{decayW} and in the following $C$ denotes various positive constant which may vary from line to line but are independent of $x$.
We also have on $\left[0,T/2\right]$
\begin{equation}\label{evolution-1}
|\dot x_1|\leq C W^{\frac{\nu}{\alpha}}.
\end{equation}
On $\left(T/2,T\right)$ the closed-loop system is
\begin{equation}
\dot{x}_1  =   - k_1 \{ x_1 \}^{1+\kappa_1},\;
\dot{x}_2  =  - k_1 \{ x_1 \}^{1+\kappa_1} x_3,\;
\dot x_3  =   0.
\label{cl2_unicycle-simp1}
\end{equation}
In particular, if $V(x_1):=x_1^2$,
\begin{equation}\label{dotVx1}
\dot V\leq -\frac{1}{C}V^{\frac{2+\kappa_1}{2}}.
\end{equation}
We also have on $\left(T/2,T\right)$
\begin{equation}\label{deuxiemepartie2et3}
|\dot x_2|+ |\dot x_3| \leq C(|x_2|+|x_3|)V^{\frac{1+\kappa_1}{2}}.
\end{equation}
The conclusion of Theorem~\ref{thm:unicycle_LFTS_1} easily follows
from \eqref{decayW}, \eqref{evolution-1}, \eqref{dotVx1}, and \eqref{deuxiemepartie2et3}. (Let us recall that $\kappa_1<0$ and $\kappa_2<0$.)

\hfill $\Box$

Let us now return to the initial control system \eqref{unicyclexcossin}. Let us deduce from Theorem~\ref{thm:unicycle_LFTS_1} the following theorem.

\begin{theorem}\label{thm:vrai-unicycle_LFTS_1}
Let $T>0$, $\nu\in(0,1)$, $-1<\kappa_1<0$, $-(1-\nu)/<\kappa_2<0$, $k_1>0$, and $k_2>0$. Then there exists a large enough $k_3>0$ such that, for the feedback law  $u=(u_1,u_2)\in C_{\text{t-p}} (\mathbb{R} \times \mathbb{R}^3; \mathbb{R}^2)$ defined by \eqref{u1_FTS_unicycle_1} to \eqref{u1u2_unicycle_en0_fast4} and for $f(x,u)= (u_1\cos(x_3),u_1\sin(x_3), u_2)^{\top}$,  \eqref{unulen0}, \eqref{uperiodic}, and \eqref{0toujours} hold  and there exists $\epsilon >0$ such that \eqref{epsilontoujours} holds too.
 In particular \eqref{unicyclexcossin} is  locally stabilizable in small time by means of explicit time-piecewise continuous % OK pour \modifsjm{stationary} 
 stationary periodic feedback laws.
\end{theorem}
\noindent
{\bf Proof.} The proof is a direct consequence of Corollary \ref{corcascaded} combined with the proof given above for Theorem \ref{thm:unicycle_LFTS_1}. 
% OK pour les modifs de JM \modifsjm{ 
	One can also proceed in a more direct way. Indeed simple computations and estimates show that  \eqref{decayW}, \eqref{evolution-1}, \eqref{dotVx1}, and \eqref{deuxiemepartie2et3}  still hold for \eqref{unicyclexcossin} at least if $|x_1|+|x_2|+|x_3|$ is small enough.
%A r\'{e}diger. Utiliser des arguments du type approximation homog\`{e}ne pour v\'{e}rifier que si $|x(0)|$ est tr\`{e}s petit $x_1(T/2)=0$, $x_2(T/2)=0$ et $|x_3(T/2)|$ est petit, puis que si $x_1(T/2)=0$, $x_2(T/2)=0$ et $|x_3(T/2)|$ est petit alors $x_1(t)=x_2(t)=0$ pour $t\in [T/2,T]$ et $x_3(T)=0$. Peut-\^{e}tre mettre \`{a} ce sujet un corollary du Theorem \ref{thm2} apr\`{e}s le Corollary \ref{cor2})  Insister aussi sur le fait que si $x(t)=0$ alors $x(t')=0$ pour tout $t'\geq t$

\hfill $\Box$

%\begin{figure}[!h]
%	\begin{center}
%		\includegraphics[width=0.9\textwidth]{figure_unicycle}
%	\end{center}
%	\caption{Simulation of (\ref{unicyclexcossin}) using control (\ref{u1u2_unicycle_en0_fast4}) with $T=2, \nu=0.5; \kappa_1=-0.5, \kappa_2=-0.1, l= 3, k_1=10, k_2=k_3=5.$}
%	\label{fig:simu1}
%\end{figure}

\section{Small-time stabilization of the slider}
\label{sec-slider}
Let us define $(x_1,x_2,x_3,x_4,x_5,x_6)\tr \in \R^6$ and  $(u_1,u_2)\tr\in\R^2$ by
\begin{equation}\label{defxu-slider}
x_1=:x,\; x_2:=\dot x, \; x_3:=y,\; x_4:=\dot y ,\; x_5:=\psi, \; x_6:=\dot x_5, \;
u_1:=\frac{\tau_1}{m},\; u_2:=\frac{\tau_2}{I}.
\end{equation}
Then the slider control system \eqref{eq:sousactionneinertiel} becomes the control system

\begin{equation}
\label{sliderxcossin}
\dot{x}_1  =   x_2,\;
\dot{x}_2 =   u_1\cos(x_5),\;
\dot{x}_3  =   x_4,\;
\dot{x}_4  =   u_1\sin(x_5),\;
\dot{x}_5  =   x_6,\;
\dot{x}_6  =   u_2,
\end{equation}
where the state is $x=(x_1,x_2,x_3,x_4,x_5,x_6)\tr \in \R^6$ and
the control is $(u_1,u_2)\tr \in\R^2$. Similarly to the unicycle, one first observes that if $x_3 (t_0) = x_4 (t_0) = x_5 (t_0) = x_6 (t_0) = 0$  and if, after time $t_0$, $u_2$ is equal to $0$, then $x_3 (t) = x_4 (t) = x_5 (t) = x_6 (t) = 0$. The idea is then to use
a first stationary feedback law to steer $(x_3,x_4,x_5,x_6)\tr$ to $0\in \R^4$ and then take $u_2=0$ and use $u_1$ to steer $(x_1,x_2)\tr \in \R$ to $0\in \R^2$.

As for the unicycle, let us consider the following quadratic approximation of \eqref{sliderxcossin} (see, once more, \cite{2017-Coron-Rivas-SICON})
\begin{equation}
\label{slider-simp1}
\dot{x}_1= x_2, \;
\dot{x}_2 =   u_1, \;
\dot{x}_3 =   x_4, \;
\dot{x}_4  =   u_1 x_5, \;
\dot{x}_5  = x_6,\;
\dot{x}_6 =  u_2,
\end{equation}
where the state is still $x=(x_1,x_2,x_3,x_4,x_5,x_6)\tr \in \R^6$ and
the control is still $(u_1,u_2)\tr \in\R^2$.

As for \eqref{sliderxcossin}, we have for \eqref{slider-simp1} that  if $x_3 (t_0) = x_4 (t_0) = x_5 (t_0) = x_6 (t_0) = 0$  and if, after time $t_0$, $u_2$ is equal to $0$, then $x_3 (t) = x_4 (t) = x_5 (t) = x_6 (t) = 0$. Let us therefore  choose, according to Theorem \ref{thmcascaded}, the controls $u_1$ and $u_2$ in a first step, so that the dynamics of the subsystem $(x_3, x_4, x_5, x_6 )^\top$ is such that we obtain $x_3 (t) = x_4 (t) = x_5 (t) = x_6 (t) = 0$ for all $t \geq T/2$ provided that $x_3 (0)$, $x_4 (0)$, $x_5 (0)$ and $x_6 (0)$ are small enough. Then we impose $u_2 = 0$ for $t \geq T/2$, $x_3$, $x_4$, $x_5$ and $x_6$ will remain zero. To stabilize in small time the double integrator made of the state variables $x_1$ and $x_2$, we can choose for $t > T/2$, $u_1$ as in (\ref{FLDI}).

% OK avec les modifs de JM \modifsjm{ 
Let $k_1>0$, $k_2>0$, $k_3>0$, $k_4>0$, $k_5>0$, $k_6>0$, $\kappa_1\in (-1/2,0)$,
 %$\kappa_2 \in (-\infty,0)$, 
 and $\mu \in (0,1)$. We assume that
\begin{equation}
\label{assumption-kappa}
 -\frac{1-\mu}{2(2-\mu)}<\kappa_2 <0,
\end{equation}
%\Wilfrid{On avait $\kappa_2 >-\frac{(1-\mu)^2}{(2-\mu)(3-2\mu)}$ dans la V8 mais en fait d'après le theorem 3.7 il faut $\kappa_2 >-\frac{(1-\nu)}{2}$ avec ici $\nu=\frac{\mu r_5}{r_3}=\mu (1-\mu)(1+2\kappa_2)$ ce qui donne $\kappa_2 >-\frac{1}{2},$ en regardant les autres conditions $r_7=(1-\mu)(1+2\kappa_2)+2\kappa_2>0$ cela donne $  \kappa_2 >-\frac{1-\mu}{2(2-\mu)} $ et dans ce cas on a bien $r_4,r_5,r_6,r_7>0$. et j'ai supprimé $\kappa_2 \in (-\infty,0)$ qui est maintenant explicite dans cette condition} and define
\begin{gather}
\label{derri}
r_1:=1,\, r_2:=1+\kappa_1, \, r_3:=1, \, r_4:=r_3+\kappa_2,\, r_5:=(r_4+\kappa_2)(1-\mu),\, r_6:= r_5+\kappa_2,
\, r_7:=r_6+\kappa_2,
\\
\label{defbarx4}
\bar x_4:=-k_3\pui{x_3}{\frac{r_4}{r_3}},
\\
\label{defbarx5}
\bar x_5:=-k_4\pui{\pui{x_4}{\frac{r_3}{r_4}}-\pui{\bar x_4}{\frac{r_3}{r_4}}}{\frac{r_5}{r_3}},
\\
\label{defbarx6}
\bar x_6:=-k_5\pui{\pui{x_5}{\frac{r_3}{r_5}}-\pui{\bar x_5}{\frac{r_3}{r_5}}}{\frac{r_6}{r_3}}.
\end{gather}
Note that \eqref{assumption-kappa} insures that $r_4>0$, $r_5>0, r_6>0$, and $r_7>0$. Let us then prove the following theorem, ensuring the small-time stabilization of the origin for  (\ref{slider-simp1}) with explicit time-varying feedback laws.

\begin{theorem}\label{thm:slider_LFTS}
	Let $T>0$, $k_1>0$, $k_2>0$, $k_3>0$, $k_4>0$, $\kappa_1\in (-1/2,0)$,
 %$\kappa_2 \in (-\infty,0)$, 
 $\mu \in (0,1)$ and $\kappa_2<0$ be given so that \eqref{assumption-kappa} is satisfied. We assume  that
\begin{equation}\label{condition-k4}
k_4>g\left(\frac{1}{1+\kappa_2},\kappa_2,k_3\right)=2^{-2\kappa_2}k_3^{\frac{1}{1+\kappa_2}},
\end{equation}
%\Wilfrid{Avant on avait $2^{-2\kappa_2(1+\kappa_2)}k_3^{\frac{1}{1+\kappa_2}}$ en refaisant les calculs on a la formule ci-dessus}
(let us recall that $g$ is defined in \eqref{eq:condition_gain}). Then there existe $K_5>0$ such that, for every $k_5 \geq K_5$,
 there exists $K_6$ such that, for every $k_6\geq K_6$,  the following feedback law  $u=(u_1,u_2)\in C_{\text{t-p}} (\mathbb{R} \times \mathbb{R}^6; \mathbb{R}^2)$:
\begin{gather}
\label{def_u-1-part-1}
		u_1 = \left|\pui{x_4}{\frac{r_3}{r_4}}-\pui{\bar x_4}{\frac{r_3}{r_4}}\right|^{\frac{\mu r_5}{r_3}} \text{ in } \left[0,T/2\right]\times \mathbb{R}^6,
\\
\label{def_u-2-part-1}
		u_2 =  -k_6 \pui{\pui{x_6}{\frac{r_3}{r_6}} - \pui{\bar x_6}{\frac{r_3}{r_6}}}{\frac{r_6+\kappa_2}{r_3} }\text{ in } \left[0,T/2\right]\times \mathbb{R}^6, \;
\\
\label{def_u-1-part-2}
u_1 = -k_{1}\pui{x_1}{\frac{r_2+\kappa_1}{r_1}}-k_2\pui{x_2}{\frac{r_2+\kappa_1}{r_2}} \,  \text{ in } \left(T/2,T\right)\times \mathbb{R}^6,
\\
\label{def_u-2-part-2}u_2 = 0 \text{ in } \left(T/2,T\right)\times \mathbb{R}^6,
\\
\label{u1u2-eriodic}
u(t+T,x)=u(t,x),\; \forall t \in \mathbb{R},\; \forall x \in \mathbb{R}^6,
\end{gather}
is such that, for the homogeneous approximation (\ref{slider-simp1}) and for the slider control system \eqref{sliderxcossin},  \eqref{unulen0}, \eqref{uperiodic}, and \eqref{0toujours} hold,  and there exists $\epsilon >0$ such that, for
\begin{gather}
\label{f-slider-homogene}
f(x,u)= (x_2,\;u_1,\; x_4,\; u_1 x_5,\;x_6 ,\; u_2)^{\top}
\end{gather}
and for
\begin{gather}
\label{f-slider}
f(x,u)= (x_2,\;u_1\cos(x_5),\; x_4,\; u_1 \sin(x_5),\;x_6 ,\; u_2)^{\top},
\end{gather}
one has
\begin{gather}
		(\dot x = f(x,u(t,x)) \;\text{and} \; |x(s)|\leq \epsilon) \Rightarrow (x(\tau) = 0), \; \forall \tau \geq s+2T), \; \forall s \in \mathbb{R}.
\end{gather}
In particular \eqref{sliderxcossin} is  locally stabilizable in small time by means of explicit time-piecewise continuous % OK avec la modif \modifsjm{stationary} 
stationary periodic feedback laws.
% OK avec les \commentsjm{Les feedbacks consid\'{e}r\'{e}s dans la version pr\'{e}c\'{e}dentes sont tr\`{e}s diff\'{e}rents. Mais je ne crois pas qu'avec les th\'{e}or\`{e}mes montr\'{e}s ici ou pr\'{e}c\'{e}demment connus on puisse montrer que les anciens feedbacks conviennent.}
\end{theorem}

\noindent
 \textbf{Proof of Theorem~\ref{thm:slider_LFTS}.} The proof is an application of Theorem \ref{thmcascaded}, Theorem~\ref{th-robustness-stability}, Theorem~\ref{thmG1} and Theorem~\ref{thmDI}. Let us first deal with (\ref{slider-simp1}).
% OK avec les \modifsjm{
 On $\left[0,T/2\right]$, the quadartic approximation (\ref{slider-simp1}) reads as
\begin{equation}
\label{x-3-x-4-x-5-x-6}
\left\{
\begin{array}{l}
\dot{x}_3 =   x_4,\\
\dot{x}_4  =   \left|\pui{x_4}{\frac{r_3}{r_4}}-\pui{\bar x_4}{\frac{r_3}{r_4}}\right|^{\frac{\mu r_5}{r_3}} x_5,\\
\dot{x}_5  = x_6,\\
\dot{x}_6 =  -k_6 \pui{\pui{x_6}{\frac{r_5}{r_6}} - \pui{\bar x_6}{\frac{r_5}{r_6}}}{\frac{r_6+\kappa_2}{r_5} }.\\
\end{array}
\right.
\end{equation}
Let us consider the control system
\begin{equation}
\dot{x}_3 =   x_4,\, \dot{x}_4 =\left|\pui{x_4}{\frac{r_3}{r_4}}-\pui{\bar x_4}{\frac{r_3}{r_4}}\right|^{\frac{\mu r_5}{r_3}}x_5,
\end{equation}
where the state is $(x_3,x_4)\tr \in \R^2$ and the control is $x_5\in \R$. It follows from Theorem~\ref{thmDI} and \eqref{condition-k4}  that $0\in \R^2$ is small-time stable for this control system if one uses the feedback law $x_5=\bar x_5$.  Let us point out that if
\begin{equation}\label{def-f}
f(x_3,x_4,x_5)=\left(x_4, \left|\pui{x_4}{\frac{r_3}{r_4}}-\pui{\bar x_4}{\frac{r_3}{r_4}}\right|^{\frac{\mu r_5}{r_3}} x_5\right)\tr,
\end{equation}
one has
\begin{gather}
\label{hom_f1}
f_1(\lambda^{r_3}x_3,\lambda^{r_4}x_4,\lambda^{r_5}x_5)
=\lambda^{r_3+\kappa_2}f_1(x_3,x_4,x_5),
\\
\label{hom_f2}
f_2(\lambda^{r_3}x_3,\lambda^{r_4}x_4,\lambda^{r_5}x_5) =\lambda^{r_4+\kappa_2}f_2(x_3,x_4,x_5).
\end{gather}
Moreover, for $l :=r_3/r_6$, one has, using \eqref{assumption-kappa} and \eqref{derri},
\begin{gather}
l+1>\frac{r_3}{r_6}>\frac{r_4}{r_6}>\frac{r_5}{r_6}>0,
\\
\label{hom-x5}
\bar x_5(\lambda^{r_3}x_3,\lambda^{r_4}x_4)=\lambda^{r_5}\bar x_5(x_3,x_4),
\\
\pui{\bar x_5}{l} \text{ is of class }  C^1.
\end{gather}
Then, by Theorem~\ref{thmcascaded}, at least if $k_5>0$ is large enough, $0\in \R^3$ is small-time stable
for
\begin{equation}
\label{x-3-x-4-x-5-x-6-closed}
\left\{
\begin{array}{l}
\dot{x}_3 =   x_4,\\
\dot{x}_4  =   \left|\pui{x_4}{\frac{r_3}{r_4}}-\pui{\bar x_4}{\frac{r_3}{r_4}}\right|^{\frac{\mu r_5}{r_3}} x_5,\\
\dot{x}_5  = -k_5 \pui{\pui{x_5}{\frac{r_3}{r_5}} - \pui{\bar x_5}{\frac{r_3}{r_5}}}{\frac{r_6}{r_3} }=\bar x_6.
\end{array}
\right.
\end{equation}
Note that $\pui{\bar x_6}{l^\prime}$ is of class $C^1$ for $l^\prime:=r_3/r_7$.
One applies once more  Theorem~\ref{thmcascaded} (see also Remark~\ref{rem-Morin-Samson}). Using once more  \eqref{assumption-kappa} and \eqref{derri}, one gets
\begin{gather}
r_7=r_6+\kappa_2>0,
\\
l+1=\frac{r_3}{r_6}+1 >\frac{r_3}{r_6+\kappa_2}>\frac{r_4}{r_6+\kappa_2}>\frac{r_5}{r_6+\kappa_2}
 >\frac{r_6}{r_6+\kappa_2}>0.
\end{gather}
One gets that,
under the assumptions of Theorem~\ref{thm:slider_LFTS}, $0\in \R^4$ is small-time stable  for \eqref{x-3-x-4-x-5-x-6-closed}. 
By Theorem~\ref{th-robustness-stability} one gets that, under the same assumptions, $0\in \R^4$ is small-time stable  for
\begin{equation}
\label{x-3-x-4-x-5-x-6-closed-slider}
\left\{
\begin{array}{l}
\dot{x}_3 =   x_4,\\
\dot{x}_4  =   \left|\pui{x_4}{\frac{r_3}{r_4}}-\pui{\bar x_4}{\frac{r_3}{r_4}}\right|^{\frac{\mu r_5}{r_3}} \sin(x_5),\\
\dot{x}_5  = -k_5 \pui{\pui{x_5}{\frac{r_3}{r_5}} - \pui{\bar x_5}{\frac{r_3}{r_5}}}{\frac{r_6}{r_3} }=\bar x_6.
\end{array}
\right.
\end{equation}
%}
In particular if $x:[0,T/2] \to \R^6$ is a solution of the closed-loop systems $\dot x =f(x,u(t,x))$ where $f$ is given by \eqref{f-slider-homogene} or by \eqref{f-slider}, we have $x_3(T/2)=x_4(T/2)=x_5(T/2)=x_6(T/2)=0$ if
$|x_3(0)|+|x_4(0)|+|x_5(0)|+|x_6(0)|$ is small enough.

%OK pour le smodifications \modifsjm{
On $[T/2, T]$, one has, if $x_5(T/2)=x_6(T/2)=0$ for $f$ given by \eqref{f-slider} and for every $(x_5(T/2),x_6(T/2)) $ for $f$ given by \eqref{f-slider-homogene},
\begin{equation}\label{slider-surT2T}
\dot x_1=x_2, \,  \dot x_2= -k_{1}\pui{x_1}{\frac{r_2+\kappa_1}{r_1}}-k_2\pui{x_2}{\frac{r_2+\kappa_1}{r_2}}.
\end{equation}
Note that, by Theorem~\ref{thmG1}, $0\in \R^2$ is small-time stable for \eqref{slider-surT2T}. 
Theorem~\ref{thm:slider_LFTS} follows from the above arguments.
%}
% OK avec la modif \modifsjm{
\begin{remark}
In Theorem~\ref{thm:slider_LFTS}, for $t\in (T/2,T)$ one can replace $u_1$ defined in \eqref{def_u-1-part-2} by $u_1=u$ where $u$ is defined in \eqref{FLDI} provided that $l$, $\kappa$, $k_1$, and $k_2$ satisfy the assumptions of Theorem~\ref{thmDI}. Indeed it suffices to then use in the above proof Theorem~\ref{thmDI} instead of Theorem~\ref{thmG1}. 
\end{remark}
%}

 \hfill $\Box$

\section*{Conclusion}
To conclude, new explicit feedback laws, based on cascaded and desingularization techniques, and homogeneity properties, have been proposed leading to small-time stabilization of two classical examples of nonholonomic or underactuated mechanical systems: the unicycle robot and the slider, which are known to be not asymptotically stabilizable by means of continuous feedback laws. The proposed resulting feedback laws have been first established considering homogeneous approximations of the mechanical systems  and then successfully extended to the original ones.

\appendix
\section{Some useful inequalities}
In this appendix we prove two lemmas which are used in previous  sections.

\begin{lemma}\label{lem:ineq1} Let $l >1$, then $\forall x,y\in\R : |x-y|^2\leq 2^{2(l-1)/l}|\pui{x}{l}-\pui{y}{l}|^{2/l}$.
\end{lemma}
\noindent
\textbf{Proof.}
This inequality is true for $x=y$. Assume $(x-y)\neq 0$. For any given $l>1$, let us define for $x,y \in \R$ ($x\neq y$)  $g_l(x,y) = \frac{\vert \pui{x}{l}-\pui{y}{l}\vert^{2/l}}{(x-y)^2}$. We have $g_l(x,y)  = h_l\left(\frac{y}{x-y}\right)$ with $h_l : z\in \R \mapsto \vert \pui{1+z}{l}-\pui{z}{l}\vert^{2/l}$ which is bounded below by $2^{2(1-l)/l}$ which gives the desired inequality.
\hfill $\Box$

\begin{lemma}\label{lem:ineq2} Let $a_0>0$, $a_1>0$, $a_2>0$, and $\beta>\gamma>1$, then
\begin{equation}\label{psi<0}
\psi(z)=-a_0-a_1z^{\beta}+a_2z^{\gamma} <0, \, \forall z \in [0,+\infty),
\end{equation}
if and only if
\begin{equation}\label{cnspsi>0}
a_0 a_1^{\frac{\gamma}{\beta-\gamma}}> a_2^{\frac{\beta}{\beta-\gamma}}
\left(\left(\frac{\gamma}{\beta}\right)^{\frac{\gamma}{\beta-\gamma}}
-\left(\frac{\gamma}{\beta}\right)^{\frac{\beta}{\beta-\gamma}}\right).
\end{equation}
\end{lemma}
\begin{remark}
\label{rem.simplecondition}
Let $a_0>0$, $a_1>0$, $a_2>0$, $\beta \in \R$, and $\gamma \in \R\setminus\{0\}$. Then
\begin{equation}\label{simplecondition}
 a_1> \frac{\gamma}{\beta}\frac{a_2^\frac{\beta}{\gamma}}{a_0^{\frac{\beta-\gamma}{\gamma}}}
\end{equation}
implies \eqref{cnspsi>0}.
\end{remark}
\noindent
\textbf{Proof of Lemma~\ref{lem:ineq2}.}
$\psi^\prime(z)=-a_1\beta z^{\beta-1}+a_2\gamma z^{\gamma-1}$ which zeros at $0$ and at $z_0=\left(\frac{a_2\gamma}{a_1\beta}\right)^{\frac{1}{\beta-\gamma}}$ thus $\psi$ on $[0,+\infty[$ has a maximum at $z_0$:  $\psi(z_0)=-a_0-a_1\left(\frac{a_2\gamma}{a_1\beta}\right)^{\frac{\beta}{\beta-\gamma}}+a_2\left(\frac{a_2\gamma}{a_1\beta}\right)^{\frac{\gamma}{\beta-\gamma}}$  which is negative  if and only if \eqref{cnspsi>0} holds.
\hfill $\Box$

\bibliographystyle{plain}

\end{document}